\documentclass{amsart}
\usepackage{amssymb,amsmath}
\usepackage[american,french]{babel}
\usepackage{amsopn}
\usepackage{fancyhdr}

\def\Titre{Exponential mixing for the $3D$ stochastic Navier--Stokes equations}
\title{\Titre}
\author{Cyril ODASSO
\\
 \\  Ecole Normale Sup\'erieure de Cachan, antenne de Bretagne,\\ Avenue Robert Schuman,
 Campus de Ker Lann, 35170 Bruz (FRANCE). \\ and
\\ IRMAR,  UMR 6625 du CNRS, Campus de Beaulieu,  35042 Rennes cedex (FRANCE)\\}

\newtheorem{Theorem}{Theorem}[section]
\newtheorem{Definition}[Theorem]{Definition}
\newtheorem{Proposition}[Theorem]{Proposition}
\newtheorem{Lemma}[Theorem]{Lemma}
\newtheorem{Corollary}[Theorem]{Corollary}
\newtheorem{Remark}[Theorem]{Remark}

\newtheorem{Hypothesis}[Theorem]{Hypothesis}
\makeatletter
\newif\ifmsbmloaded@

\@addtoreset{equation}{section}

\makeatother

\def\R{\mathbb R}
\def\N{\mathbb N}

\def\E{\mathbb E}
\def\P{\mathbb P}
\def\Pcal{\mathcal{P}}

\def\H{\mathbb H}

\def\Dr{\mathcal D}
\def\Br{\mathcal B}

\def\F{\mathcal F}

\def\H{\Bbb H}

\newcommand{\BLANC}[1]{   }

\newcommand{\abs}[1]{\left\vert#1\right\vert}
\newcommand{\norm}[1]{\left\Vert#1\right\Vert}
\newcommand{\eps}{\varepsilon}
\newcommand{\sig}{\sigma}

\def \Espace{\renewcommand{\arraystretch}{1.7} }

\begin{document}
\setcounter{page}{205}

\selectlanguage{american}

\maketitle

\pagestyle{fancy}

\noindent\textbf{Abstract}:
We study  the Navier-Stokes equations in dimension $3$ (NS3D) driven by a noise which is white in time.
We establish that if the noise is at same time sufficiently smooth and non degenerate in space,
then the weak solutions converge exponentially fast to equilibrium.

We use a coupling method. The arguments used in dimension two do
not apply since, as is well known, uniqueness is an open problem
for NS3D. New ideas are introduced. Note however that
 many simplifications appears since we work with non degenerate noises.

\

\noindent {\bf Key words}: Stochastic three-dimensional Navier-Stokes equations,
Markov transition semi-group, invariant measure, ergodicity, coupling method,
exponential mixing, galerkin approximation.

\section*{Introduction}

We are concerned with the stochastic Navier--Stokes equations on a three dimensional bounded domain (NS3D)
 with Dirichlet boundary conditions.
These equations describe the time evolution of an incompressible fluid subjected to a determinist and a random exterior
force and are given by
\begin{equation}\label{EqIntroNS}
\Espace
\left\{
\begin{array}{l}
dX+\nu(-\Delta)X \, dt+(X,\nabla)X\,dt+\nabla p \, dt =  \phi(X)dW+ f \, dt,\\
\begin{array}{rclll}
        \left(\textrm{div } X\right)(t,\xi) &=& 0, &\textrm{ for } \xi\in D,&t>0,\\
        X(t,\xi) &=& 0,&\textrm{ for } \,\xi\in\partial D,&t>0,\\
        X(0,\xi)&=& x_0(\xi), &\textrm{ for } \xi\in D.&
\end{array}
\end{array}
\right.
\end{equation}
Here $D$ is an open bounded domain of $\R^3$ with smooth boundary $\partial D$ or $D=(0,1)^3$.
We have denoted by $X$ the velocity, by $p$ the pressure and by $\nu$ the viscosity.
 The  external force field acting on the fluid  is the sum of a
random force field of white noise type $\phi(X)dW$ and a determinist one $f \, dt$.

In the deterministic case ($\phi=0$), there exists a global weak
solution (in the PDE sense) of \eqref{EqIntroNS} when $x_0$ is
square integrable,
 but uniqueness of such solution is not known. On another hand, there exists a unique local strong solution
 when $x_0$ is smooth, but global existence is an open problem (see \cite{temam-survey} for a survey on these
questions).

In the stochastic case, there exists a global weak solution of the
martingale problem, but pathwise uniqueness or uniqueness in law
remain open problems (see \cite{flandoli-cetraro} for a survey on
the stochastic case).

The main result of the present article is to establish that, if
$\phi$ is at the same time sufficiently smooth and non degenerate,
then the solutions converge exponentially fast to equilibrium.
More precisely, given a solution, there exists a stationary
solution (which might depends on the given solution), such that
the total variation distance between the laws of the given
solution and of the stationary solution converges to zero
exponentially fast.

Due to the lack of uniqueness, it is not straightforward to define
a Markov evolution associated to \eqref{EqIntroNS}. Some recent
progress have been obtained in this direction. In
\cite{DebusscheNS3D}, \cite{ODASSOMarkov},  under conditions on
$\phi$ and $f$ very similar to ours, it is shown that every
solution of \eqref{EqIntroNS} limit of Galerkin approximations
verify the weak Markov property. Uniqueness in law is not known
but we think that this result  is a step in this direction. Our
result combined with this result implies that the transition
semi-group constructed in \cite{DebusscheNS3D} is exponentially
mixing.

Note also that recently, a Markov selection argument has allowed
the construction of a Markov evolution in \cite{FlandoliRomito2}.
Our result does not directly apply since we only consider
solutions which are limit of Galerkin approximations. However,
suitable modifications of our proof might imply that under
suitable assumptions on the noise, the Markov semi-group
constructed in \cite{FlandoliRomito2} is also exponentially
mixing.

Our proof relies on coupling arguments. These have been introduced
recently in the context of stochastic partial differential
equations by several authors (see \cite{H}, \cite{K}, \cite{KS},
\cite{KS2}, \cite{KS3},  \cite{Matt}, \cite{ODASSO1},
\cite{ODASSO3} and \cite{S}). The aim was to prove exponential
mixing for degenerate noise. It was previously observed that the
degeneracy of the noise on some subspace could be compensated by
dissipativity arguments \cite{BKL}, \cite{EMS}, \cite{KS00}. More
recently, highly degenerate noise noises have been considered in
\cite{HM04}, \cite{MattPar}.

In all these articles, global well posedness of the stochastic
equation is strongly used in many places of the proof. As already
mentioned, this is not the case for the three dimensional
Navier-Stokes equations considered here. Thus substantial changes
in the proof have to be introduced. However, we require that the
noise is sufficiently non degenerate and many difficulties of the
above mentioned articles disappear.

The main idea is that coupling of solutions can be achieved for
initial data which are small in a sufficiently smooth norm.
 A coupling satisfying good
properties is constructed thanks to the Bismut-Elworthy-Li
formula. Another important ingredient in our proof is that any
weak solution enters a small ball in the smooth norm and that the
time of entering in this ball admits an exponential moment. We
overcome the lack of uniqueness of solutions by working with
Galerkin approximations. We prove exponential mixing for these
with constants which are controlled uniformly. Taking the limit,
we obtain our result for solutions which are limit of Galerkin
approximations.

\pagebreak
\section{Preliminaries and main result}

\subsection{Weak solutions}
\

Here $\mathcal L(K_1;K_2)$ (resp $\mathcal L_2(K_1;K_2)$) denotes
the space of bounded (resp Hilbert-Schmidt) linear operators from
the Hilbert space $K_1$ to $K_2$.

 We denote by $\abs\cdot$ and $(\cdot,\cdot)$ the norm and the inner product of $L^2(D;\R^3)$ and by $\abs\cdot_p$ the norm of $L^p(D;\R^3)$.
Recall now the definition of the Sobolev spaces $H^p(D;\R^3)$ for
$p\in \N$
$$
\Espace
\left\{
\begin{array}{l}
H^p(D;\R^3)=\left\{X\in L^2(D;\R^3)\,\left|\,\partial_\alpha X \in L^2(D;\R^3)\,\textrm{ for } \abs{\alpha}\leq p\right.\right\},\\
\abs{X}^2_{H^p}=\sum_{\abs{\alpha}\leq p}\abs{\partial_\alpha X}^2.
\end{array}
\right.
$$
It is well known that $(H^p(D;\R^3),\abs{\cdot}_{H^p})$ is a Hilbert space. The Sobolev space $H^1_0(D;\R^3)$ is
the closure of the space of smooth functions on $D$ with compact support  by $\abs\cdot_{H^1}$. Setting
$
\norm{X}=\abs{\nabla X},
$
we obtain that $\norm\cdot$ and $\abs\cdot_{H^1}$ are two equivalent norms on $H^1_0(D;\R^3)$ and that
$(H^1_0(D;\R^3),\norm\cdot)$ is a Hilbert space.

 Let $H$ and $ V $ be the closure of the space of
smooth functions on $D$ with compact support and free divergence
for the norm $\abs \cdot$ and $\norm\cdot$, respectively.

\noindent Let $\pi$ be the orthogonal projection in $L^2(D;\R^3)$ onto the space $H$. We set
$$
A = \pi\left(-\Delta\right), \,  D(A) =  V  \cap H^2(D;\R^3), \,
B(u,v)=\pi \left((u,\nabla)v\right)\,\textrm{ and
}\,B(u)=B(u,u).$$ Let us recall the following useful identities
$$
\Espace\left\{ \begin{array}{rclll}
 (B(u,v),v)&=&0, &u,\,v \in V,\\
 (B(u,v),w)&=&-(B(u,w),v),&u,\,v,\,w\in V.
 \end{array}
 \right.
$$

\noindent  As is classical, we get rid of the pressure and rewrite problem \eqref{EqIntroNS} in the form
\begin{equation}\label{EqNS}
\Espace
\left\{
\begin{array}{rcl}
dX+ \nu A X dt+B(X) dt &=& \phi(X)dW+f \, dt,\\
                             X(0)&=& x_0,
\end{array}
\right.
\end{equation}
where $W$ is a cylindrical Wiener process on $H$ and with a slight abuse of notations, we have denoted by the same symbols
 the projections of $\phi$ and $f$.

  It is well-known  that $\left(A,\Dr(A)\right)$ is
 a self-adjoint operator with discrete spectrum. See \cite{Cons}, \cite{Temam}.
We consider $(e_n)_n$ an eigenbasis of $H$ associated to the increasing sequence $(\mu_n)_n$ of eigenvalues
 of $\left(A,\Dr(A)\right)$.
It will be convenient to use  the fractionnal power $\left(A^s,\Dr(A^s)\right)$ of the operator
$(A,\Dr(A))$ for $s\in \R$
$$
\Espace
\left\{
\begin{array}{rcl}
\Dr(A^s)&=&\left\{X=\sum_{n=1}^\infty x_n e_n\,\left|\,\sum_{n=1}^\infty\mu_n^{2s}\abs{x_n}^2<\infty\right.\right\},\\
A^s X&=&\sum_{n=1}^\infty \mu_n^s x_n e_n\;\textrm{ where }\;X=\sum_{n=1}^\infty x_n e_n.
\end{array}
\right.
$$
We set for any $s\in \R$
$$
\norm{X}_s=\abs{A^\frac{s}{2}X},\quad
\H_s=\Dr(A^\frac{s}{2}).
$$
It is obvious that $(\H_s,\norm\cdot_s)$ is a Hilbert space, that $(\H_0,\norm\cdot_0)=(H,\abs\cdot)$ and
that $(\H_1,\norm\cdot_1)=(V,\norm\cdot)$. Moreover, recall that, thanks to the regularity theory of the Stokes operator,
$\H_s$ is a closed subspace of $H^s(D,\R^3)$ and $\norm\cdot_s$ is equivalent to the usual norm
 of $H^s(D;\R^3)$ when $D$ is an open bounded domain of $\R^3$ with smooth boundary $\partial D$.
When $D=(0,1)^3$, it remains true for $s\leq 2$.

Let us define
$$
\Espace
\left\{
\begin{array}{rcl}
\mathcal X&=&L^\infty_{\textrm{loc}}(\R^+;H)\cap L^2_{\textrm{loc}}(\R^+;V)\cap C(\R^+;\H_s),\\
\mathcal W&=&C(\R^+;\H_{-2}),\\
{\Omega_*}&=& \mathcal X\times \mathcal W,
\end{array}
\right.
$$
where $s$ is any fixed negative number. Remark that the definition of $\mathcal X$  is not depending on $s<0$.
 Let ${X_*}$ (resp ${W_*}$) be the projector
 ${\Omega_*}\to\mathcal X$ (resp ${\Omega_*}\to\mathcal W$). The space ${\Omega_*}$ is endowed with its
  borelian $\sig$-algebra
${\F^*}$ and with $\left({\F_t^*}\right)_{t\geq 0}$ the filtration
generated by $({X_*},{W_*})$.

 Recall that $W$
 is said to be a $\left(\mathcal F_t\right)_{t}$--cylindrical Wiener process on $H$ if
$W$ is $\left(\mathcal F_t\right)_{t}$--adapted, if $W(t+\cdot)-W(t)$ is independant of
$\mathcal F_t$ for any $t\geq 0$ and if $W$ is a cylindrical Wiener process on $H$. Let $E$ be a
Polish space. We denote by $P(E)$ the set of probability measure on $E$ endowed with the borelian $\sig$--algebra.

\begin{Definition}[Weak solutions]
 A probability measure $\P_{\lambda}$ on $\left({\Omega_*},{\F^*}\right)$ is said to be a weak solution of \eqref{EqNS}
 with initial law $\lambda\in P(H)$ if the three following properties hold.
\begin{itemize}
\item[i)] The law of ${X_*}(0)$ under $\P_{\lambda}$ is $\lambda$.

\item[ii)] The process ${W_*}$  is a $({\F_t^*})_{t}$--cylindrical
Wiener process on $H$ under $\P_{\lambda}$.

\item[iii)] We have $\P_{\lambda}$-almost surely
\begin{equation}\label{Def_weak}
\Espace
\begin{array}{r}
({X_*}(t),\psi)+\nu\int_0^t({X_*}(s),A\psi)ds+\int_0^t(B(X_*(s)),\psi)ds\quad\quad\quad\quad\quad\\
=({X_*}(0),\psi)+t\left(f,\psi\right)
+\int_0^t(\psi,\phi({X_*}(s))d{W_*}(s)),
\end{array}
\end{equation}
for any $t\in\R^+$ and any $\psi$ smooth mapping on $D$ with
compact support and divergence zero.
\end{itemize}
When the initial value $\lambda$ is not specified, $x_0$ is the initial value of the weak solution $\P_{x_0}$
 (i.e. $\lambda$ is equal to $\delta_{x_0}$ the Dirac mass at point $x_0$).
\end{Definition}
These solutions are weak in both probability and PDE sense. On the
one hand, these are solutions in law. Existence of solutions in
law does not imply that, given a Wiener process $W$ and an initial
condition $x_0$, there exist a solution $X$ associated to $W$ and
$x_0$. On the other hand, these solutions live in $H$ and it is
not known if they live in $\H_1$. This latter fact causes many
problems when trying to apply Ito Formula on $F({X_*}(t))$ when
$F$ is a smooth mapping. Actually, we do not know if we are
allowed to apply it.

That is the reason why we do not consider any weak solution but
only those which are limit in  distribution of solutions of
Galerkin approximations of \eqref{EqNS}. More precisely, for any
$N\in\N$, we denote by $P_N$ the eigenprojector of $A$ associated
to the first $N$ eigenvalues. Let $(\Omega,\F,\P)$ be a
probability space and $W$ be a cylindrical Wiener process on $H$
for $\P$. We consider the following approximation of \eqref{EqNS}
\begin{equation}\label{galerkin}
\Espace
\left\{
\begin{array}{rcl}
dX_N+  \nu A X_N dt+P_NB(X_N) dt &=& P_N\phi(X_N)dW+P_Nfdt,\\
                             X_N(0)&=& P_N x_0.
\end{array}
\right.
\end{equation}
In order to have existence of a weak solution, we use the
following assumption.
\begin{Hypothesis}
\label{HGal}
The mapping $\phi$
is bounded Lipschitz $H\to \mathcal L_2\left(H;\H_1 \right)$
and $f\in H$.
\end{Hypothesis}
We set
$$
B_1=\sup_{x\in H}\abs{\phi(x)}_{\mathcal L_2\left(H;\H_1
\right)}^2+\frac{\abs{f}^2}{\nu \mu_1}.
$$

It is easily shown that, given $x_0\in H$, 
 \eqref{galerkin} has a unique solution $X_N=X_N(\cdot,x_0)$. Proceeding as in
\cite{FlandoliGatarek}, we can see that the laws $(\P_{x_0}^N)_N$
of $(X_N(\cdot,x_0),W)$ are tight in a well chosen functional
 space. Then, for a subsequence $(N_k)_k$, $(X_{N_k},W)$ converges in law to $\P_{x_0}$ a weak solution
of \eqref{EqNS}. Hence we have existence of the weak solutions of \eqref{EqNS}, but uniqueness remains an open problem.

\begin{Remark}
We only consider weak solutions constructed in that way. This
allows to make some computations and to obtain many estimates.
 For instance, when trying to estimate the $L^2$-norm of ${X_*}(t)$ under a weak solution $\P_{x_0}$,
 we would like to
apply the Ito Formula on $\abs{{X_*}}^2$. This would give
$$
d\abs{{X_*}}^2+2\nu\norm{{X_*}}^2dt=2\left({X_*},\phi({X_*})d{W_*}\right)+2(f,{X_*})dt+\abs{\phi({X_*}(t))}^2_{\mathcal
L_2(H;H)}dt.
$$
Integrating and taking the expectation, we would deduce that, if $f=0$ and $\phi$ constant,
$$
\E_{x_0}\left(\abs{{X_*}(t)}^2+2\nu\int_0^t\norm{{X_*}(s)}^2dt\right)=\abs{x_0}^2+t\abs{\phi}^2_{\mathcal
L_2(H;H)}.
$$
Unfortunately, those computations are not allowed. However,
analogous computations are valid if we replace $\P_{x_0}$ by
$\P_{x_0}^N$, which yields
$$
\E\left(\abs{X_N(t)}^2+2\nu\int_0^t\norm{X_N(s)}^2dt\right)=\abs{P_N x_0}^2+t\abs{P_N\phi}^2_{\mathcal L_2(H;H)}.
$$
 Then, we take the limit and we infer
from Fatou Lemma and from the semi-continuity of $\abs\cdot$, $\norm\cdot$ in $\H_s$
that
$$
\E_{x_0}\left(\abs{{X_*}(t)}^2+2\nu\int_0^t\norm{{X_*}(s)}^2dt\right)\leq\abs{x_0}^2+t\abs{\phi}^2_{\mathcal
L_2(H;H)},
$$
provided $f=0$ and $\phi$ constant and
provided $\P_{x_0}$ is limit in  distribution of
solutions of \eqref{galerkin}.
\end{Remark}
Let $\P'$ and $Y$ be a probability measure and a random variable
on $({\Omega_*},{\F^*})$, respectively.
 The distribution $\Dr_{\P'}(Y)$ denotes the law of $Y$ under $\P'$.

A weak solution $\P_{\mu}$ with initial law $\mu$ is said to be stationary if, for any $t\geq 0$,
 $\mu$ is equal to $\Dr_{\P_{\mu}}({X_*}(t))$.

 We define
$$
\left(\Pcal_t^N\psi\right)(x_0)=\E\left(\psi(X_N(t,x_0))\right)=\E_{x_0}^N\left(\psi({X_*}(t)\right),
$$
where $\E^N_{x_0}$ is the expectation associated to $\P_{x_0}^N$.

 It is
easily shown that $X_N(\cdot,x_0)$ verifies the strong Markov
property, which obviously implies that
$(\Pcal_t^N)_{t\in \R^+}$ is a Markov transition semi-group on $P_NH$.

Ito Formula on $\abs{X_N(\cdot,x_0)}^2$ 
 gives
$$
d\abs{X_N}^2+2\nu\norm{X_N}^2dt=2\left(X_N,\phi(X_N)dW\right)+2(X_N,f)dt+\abs{P_N\phi(X_N)}^2
dt,
$$
which yields, by applying arithmetico-geometric inequality and
Hypothesis \ref{HGal},
\begin{equation}\label{Eq1.10}
d\abs{X_N}^2+\nu\norm{X_N}^2dt\leq 2\left(X_N,\phi(X_N)dW\right)+cB_1
dt.
\end{equation}
Integrating and taking the expectation, we obtain
\begin{equation}\label{Eq1.10bis}
\E\left(\abs{X_N(t)}^2\right)\leq e^{-\nu\mu_1
t}\abs{x_0}^2+\frac{c}{\nu\mu_1}B_1.
\end{equation}
Hence, applying the Krylov-Bogoliubov Criterion (see \cite{DPZ1}),
 we obtain that $(\Pcal_t^N)_t$ admits an invariant measure $\mu_N$ and that
every invariant measure has a moment of order two in $H$. Let $X_0^N$ be a random variable whose law is $\mu_N$ and which
is independent of $W$, then $X_N=X_N(\cdot,X_0^N)$ is a stationary solution of \eqref{galerkin}.
Integrating \eqref{Eq1.10}, we obtain
$$
\E\abs{X_N(t)}^2+\nu\E\int_0^t\norm{X_N(s)}^2ds\leq\E\abs{X_N(0)}^2+cB_1t.
$$
Since the law of $X_N(s)$ is $\mu_N$ for any $s\geq 0$ and since $\mu_N$ admits a moment of order $2$, it follows
\begin{equation}\label{Eq1.11}
\int_{P_N H}\norm{x}^2\,\mu_N(dx)\leq \frac{c}{\nu}B_1.
\end{equation}
Moreover the laws $(\P_{\mu_N}^N)_N$ of $(X_N(\cdot,X_0^N),W)$ are tight in a well chosen functional
 space. Then, for a subsequence $(N_k')_k$,
 $\P_{\mu_{N_k}}^{N_k}$ converges in law to $\P_{\mu}$ a weak
stationary solution of \eqref{EqNS} with initial law $\mu$ (See \cite{FlandoliGatarek} for details).
 We deduce from \eqref{Eq1.11} that
$$
 \int_{H}\norm{x}^2\,\mu(dx)\leq \frac{c}{\nu}B_1,
$$
which yields (see \cite{FlandoliRomito})
\begin{equation}\label{Eq1.12}
\P_{\mu}\left({X_*}(t)\in\H_{1}\right)=1\;\textrm{ for any
}\;t\geq 0.
\end{equation}
We do not know if $X_*(t)\in\H_1$ for all $t$ holds
$\P_\mu$--almost surely. This would probably imply strong
uniqueness $\mu$--almost surely. Remark that it is not known in
general if $\mu$ is an invariant measure because, due to the lack
of uniqueness,
 it is not known if 
 \eqref{EqNS} defines a Markov evolution.
We will see below that this is the case under suitable
assumptions.

\subsection{Exponential convergence to equilibrium}
\


In the present article, the covariance operator $\phi$ of the
noise is assumed to be at the same time sufficiently smooth and
non degenerate with bounded derivatives. More precisely, we use
the following assumption.
\begin{Hypothesis}\label{H0}
There exist $\eps>0$  and a family $(\phi_n)_n$ of continuous
mappings $H\to\R$ with continuous derivatives such that
$$
\Espace
\left\{
\begin{array}{l}
\phi(x)dW=\sum_{n=1}^\infty\phi_n(x)e_n dW_n\quad\textrm{ where }\quad W=\sum_{n=0}^\infty W_n e_n,\\
\kappa_0=\sum_{n=1}^\infty\sup_{x\in H}\abs{\phi_n(x)}^2\mu_n^{1+\eps}<\infty.
\end{array}
\right.
$$
Moreover there exists $\kappa_1$ such that for any $x$, $\eta\in\H_2$
$$
\sum_{n=1}^\infty\abs{\phi_n'(x)\cdot \eta}^2\mu_n^{2}<\kappa_1\norm{\eta}_ 2^2.
$$
For any $x\in H$ and $N\in\N$, we have $\phi_n(x)>0$ and 
\begin{equation}\label{EqH1.4}
\kappa_2=\sup_{x\in H}\abs{\phi^{-1}(x)}_{\mathcal
L\left(\H_3;H\right)}^2<\infty,
\end{equation}
where
$$
\phi(x)^{-1}\cdot h=\sum_{n=1}^\infty \phi_n(x)^{-1} h_n
e_n\quad\textrm{ for }\quad h=\sum_{n=0}^\infty h_n e_n.
$$
\end{Hypothesis}
\noindent For instance, $\phi=A^{-\frac{s}{2}}$ fulfills Hypothesis \ref{H0} provided $s\in\left(\frac{5}{2},3\right]$.

\noindent We set
$$
B_0=\kappa_0+\kappa_1+\kappa_2+\abs{f}^2.
$$
\begin{Remark}[Additive noise]
If the noise is additive, Hypothesis \ref{H0} simplifies. Indeed
in this case, we do not need to assume that $\phi$ and $A$
commute. This requires a different but simpler proof of Lemma
\ref{PropSmallH3} below.
\end{Remark}
\begin{Remark}[Large viscosity]
Another situation where we can get rid of the assumption that the
noise is diagonal is when the viscosity $\nu$ is sufficiently
large. The proof is simpler in that case.
\end{Remark}
\begin{Remark}
It is easily shown that  Hypothesis \ref{H0} and $f\in H$ imply
Hypothesis \ref{HGal}. Therefore, solutions of \eqref{galerkin}
are well-defined and, for a subsequence, they converge to weak
solution of \eqref{EqNS}.
\end{Remark}

The aim of the present article is to establish that, under
Hypothesis \ref{H0} and under a condition of smallness of
$\norm{f}_\eps$, the law of ${X_*}(t)$ under a weak solution
$\P_{x_0}$
 converges exponentially fast to equilibrium provided $\P_{x_0}$
is limit in distribution of solutions of \eqref{galerkin}.

 Before stating our main result, let us recall some definitions.
 Let $E$ be Polish space. The set of all probability measures on
 $E$ is denoted by $\mathcal P(E)$.
The set of all bounded measurable (resp uniformly continuous) maps from $E$ to $\R$ is denoted
by $B_b(E;\R)$ (resp $UC_b(E;\R)$). The total variation $\norm{\mu}_{var}$  of a
finite real measure $\lambda$ on $E$ is given by
 $$
\norm{\lambda}_{var}= \sup  \left\{\abs{\lambda(\Gamma)}\;|\; \Gamma \in \Br(E) \right\},
$$
where we denote by $\Br(E)$ the set of the Borelian subsets of $E$.

 The main result of the present article is the following. Its
 proof is given in section $4$ after several preliminary results.
\begin{Theorem}\label{Th}
Assume that Hypothesis  \ref{H0} holds. There exists ${{\delta^0}}$, $C$ and $\gamma>0$
 only depending on $\phi$, $D$, $\eps$ and $\nu$ such that, for any weak solution $\P_{\lambda}$
 with initial law  $\lambda \in \Pcal(H)$ which  is limit of solutions of \eqref{galerkin}, there exists
a weak stationary solution $\P_{\mu}$ with initial law $\mu$ such that
\begin{equation}\label{EqTh}
\norm{\Dr_{\P_{\lambda}}({X_*}(t))-\mu}_{var}\leq C e^{-\gamma
t}\left(1+\int_{H}\abs{x}^2\,\lambda(dx)\right),
\end{equation}
provided  $\norm{f }_\eps^2\leq {{\delta^0}}$ and where
$\norm\cdot_{var}$ is the total variation norm associated to the
space $\H_s$ for $s<0$.

Moreover, for a given $\P_{\lambda}$, $\mu$ is unique and $\P_{\mu}$ is limit of solutions of \eqref{galerkin}.
\end{Theorem}
 It is well known that $\norm\cdot_{var}$ is the dual norm of $\abs\cdot_\infty$
 which means that for any finite measure $\lambda'$ on $\H_s$ for $s<0$
$$
\norm{\lambda'}_{var}=\sup_{\abs{g}_\infty\leq 1}
\abs{\int_{\H_s} g(x)\,\lambda'(dx)},
$$
where the supremum is taken over $g\in UC_b(\H_s)$ which verifies $\abs{g}_\infty\leq 1$.
Hence \eqref{EqTh} is equivalent to
\begin{equation}\label{EqTh2}
\abs{\E_{\lambda}\left(g({X_*}(t))\right)-\int_{H}g(x)\,\mu(dx)}\leq
 C\abs{g}_\infty\left(1+\int_{H}\abs{x}^2\,\lambda(dx)\right),
\end{equation}
for any $g\in UC_b(\H_s)$.
\begin{Remark}[Topology associated to the total variation norm]
Remark that if $\lambda'$ is a finite measure of $\H_{s_0}$, then
 the value of the total variation norm of $\lambda'$ associated to the
space $\H_s$ is not depending of the value of $s\leq s_0$.

Hence, since $\Dr_{\P_{\lambda}}({X_*}(t))$ is a probability
measure on $H$ then \eqref{EqTh} (resp \eqref{EqTh2}) remains true
when $\norm\cdot_{var}$ is the total variation norm associated to
the space $H$ (resp for any $g\in B_b(H; \R)$).

Moreover, we see below that, under suitable assumptions, if
$\lambda$ is a probability measure on $\H_2$, then
$\Dr_{\P_{\lambda}}({X_*}(t))$ is still a probability measure on
$\H_2$. It follows that \eqref{EqTh} (resp \eqref{EqTh2}) remains
true when $\norm\cdot_{var}$ is  associated to $\H_2$ (resp for
any $g\in B_b(\H_2; \R)$).
\end{Remark}

\BLANC{ We deduce the following result.
\begin{Corollary}[Regularization of the solutions]\label{Cor1}
Assume that Hypothesis \ref{H0} holds. There exist
 ${{\delta^0}}={{\delta^0}}(B_0,D,\eps,\nu)$, $C=C(\phi,D,\eps,\nu)>0$ and $\gamma=\gamma(\phi,D,\eps,\nu)>0$
 such that if $\norm{f }_\eps^2\leq {{\delta^0}}\,$  then, for any weak solution $\P_{\lambda}$
 with initial law  $\lambda \in \Pcal(H)$,
\begin{equation}\label{EqCor1}
\P_{\lambda}\left({X_*}(t)\not\in\H_2\right)\leq C e^{-\gamma
t}\left(1+\int_{H}\abs{x}^2\,\lambda(dx)\right),
\end{equation}
provided $\P_{\lambda}$ is limit of solutions of \eqref{galerkin}.
\end{Corollary}
The proof of this result is postponed to section $1.4$. This is a
remarkable result because ${X_*}$ living in $\H_1$ when starting
from $\H_1$ remains an open problem.
\begin{Remark}
It is well-known that Hypothesis \ref{HGal} implies that
$$
\P_{\lambda}\left({X_*}(t)\in\H_2\right)=1\quad\textrm{almost
surely in time for the Lebesgue measure},
$$
provided $\lambda\in P(H)$.

Inequality \eqref{EqCor1} of Corollary \ref{Cor1} is true for any $t\in\R^+$. Moreover, we see below that
if $\lambda\in P(\H_2)$, then
$$
\P_{\lambda}\left({X_*}(t)\in\H_2\right)=1\quad\textrm{for any
}t\in\R^+,
$$
provided $f$ and $\phi$ verifies suitable conditions.
\end{Remark}
}

Our method is not influenced by the size of the viscosity $\nu$.
Then, for simplicity in the redaction, we now assume that $\nu=1$.

\subsection{Markov evolution}
\

Here, we take into account the results of \cite{DebusscheNS3D},
\cite{ODASSOMarkov} and we rewrite Theorem \ref{Th}. This section
is not necessary in the understanding of the proof of Theorem
\ref{Th}.

Let $(N_k')_k$ be an increasing sequence of integer. In
\cite{DebusscheNS3D}, \cite{ODASSOMarkov}, it is established that
 it is possible to extract a subsequence $(N_k)_k$ of $(N_k')_k$ such
 that,  for any $x_0\in\H_2$,  $\P_{x_0}^{N_k}$
converges in distribution to a weak solution
 $\P_{x_0}$ of \eqref{EqNS} provided the following assumption holds.
\begin{Hypothesis}
\label{HCst} There exist $\eps,\delta>0$ such that the mapping
$\phi$ is bounded in $\mathcal L_2\left(H;\H_{1+\eps} \right)$.
Moreover, for any $x$,  $\textrm{ker }\phi(x)=\{0\}$ and there
exits a bounded map $\phi^{-1}: H \to\mathcal L
\left(\H_{3-\delta}; H\right)$ such that for any $x\in H$,
$$
\phi(x)\cdot \phi^{-1}(x)\cdot h=h\quad\textrm{ for any } h\in
\H_{3-\delta}.
$$
Moreover $f\in V$.
\end{Hypothesis}
 The method to
 extract $(N_k)_k$
  is based on
the investigation of the properties of the Kolmogorov equation
associated to \eqref{EqNS} perturbed by a very irregular
potential.

It follows that $(\P_{x_0})_{x_0\in \H_2}$ is a weak Markov
family, which means that for any $x_0\in \H_2$
 \begin{equation}\label{Eq1.2A}
 \P_{x_0}\left({X_*}(t)\in\H_2\right)=1\;\textrm{ for any } t\geq
0.
\end{equation}
and that, for any $t_1<\dots<t_n$, $t>0$ and any $\psi\in
B_b(\H_2;\R)$
 \begin{equation}\label{Eq1.2B}
\E_{x_0}\left(\left.\psi(X_*(t+t_n))\right|X_*(t_1),\dots,X_*(t_n)\right)=\mathcal
P_t\psi(X_*(t_n)),
\end{equation}
where
$$
\left(\Pcal_t\psi\right)(x_0)=\E_{x_0}\left(\psi({X_*}(t))\right).
$$
Note that \eqref{Eq1.2A} was known only for a stationary solution
(see \cite{FlandoliRomito}).
\begin{Remark}
Assume that  Hypothesis \ref{H0} holds. If we strengthen
\eqref{EqH1.4} into
 $$
\kappa_2=\sup_{x\in H}\abs{\phi^{-1}(x)}_{\mathcal
L\left(\H_{3-\delta};H\right)}^2<\infty,
$$
for some $\delta>0$, then Hypothesis \eqref{HCst} holds.
\end{Remark}

Hence, we immediately deduce the following corollary from Theorem
\ref{Th}.
\begin{Corollary}\label{ThDebussche1}
Assume that Hypothesis  \ref{H0} and \ref{HCst} hold. Then there
exits a unique invariant measure $\mu$
  for $\left(\Pcal_t\right)_{t\in \R^+}$ and
$C,\gamma>0$ such that for any $\lambda \in \Pcal(\H_2)$
\begin{equation}\label{EqThDebussche}
\norm{\Pcal_t^*\lambda-\mu}_{var}\leq C e^{-\gamma
t}\left(1+\int_{\H_2}\abs{x}^2\,\lambda(dx)\right),
\end{equation}
provided  $\norm{f }_\eps^2\leq {{\delta^0}}$ and where
$\norm\cdot_{var}$ is the total variation norm associated to the
space $\H_2$.
\end{Corollary}

\begin{Remark}[Uniqueness of the invariant measure $\mu$]
Assume that Hypothesis \ref{HCst} holds. Let $\P_{x_0}$ and
$\P'_{x_0}$ be two weak solutions of \eqref{EqNS} which are limit
in distribution of solutions of \eqref{galerkin}. Then we build
$(\Pcal_t)_t$ and $(\Pcal_t')_t$ as above associated to
$\P_{x_0}$ and $\P'_{x_0}$, respectively. 
 It follows that there exists $\mu$ and $\mu'$ such that \eqref{EqThDebussche} and \eqref{EqTh2} hold for
$((\Pcal_t)_t,\P_{x_0},\mu)$ and $((\Pcal_t')_t,\P'_{x_0},\mu')$.
 Although we have uniqueness of the invariant measures $\mu$ and $\mu'$ associated to $(\Pcal_t)_t$ and
 $(\Pcal_t')_t$,
 we do not know if $\mu$ and $\mu'$ are equal.
\end{Remark}

\subsection{Coupling methods}
\

 The proof of Theorem \ref{Th} is based on coupling arguments. We now recall
some basic results about coupling.
 Moreover, in order to explain the coupling method in the case of
 non degenerate noise, we briefly give the proof of exponential
 mixing for equation \eqref{galerkin}.

  Let $(\lambda_1,\lambda_2)$ be two distributions on a polish space $(E,\mathcal{E})$
 and let $(\Omega,\mathcal{F},\P)$ be a probability
 space and let $(Z_1,Z_2)$ be two random variables $(\Omega,\mathcal{F}) \to (E,\mathcal{E})$. We say that $(Z_1,Z_2)$
is a coupling of  $(\lambda_1,\lambda_2)$ if $\lambda_i=\Dr(Z_i)$ for $i=1,2$. We have denoted by
$\Dr(Z_i)$ the law of the random variable $Z_i$.

 Next result is fundamental in the coupling methods, the proof
is given for instance in the Appendix of \cite{ODASSO1}.
\begin{Lemma}\label{lem_coupling}
Let $(\lambda_1,\lambda_2)$ be two probability measures on $(E,\mathcal{E})$. Then
$$
\norm{\lambda_1-\lambda_2}_{var}= \min \P(Z_1\not = Z_2).
$$
The minimum is taken over all couplings $(Z_1,Z_2)$ of $(\lambda_1,\lambda_2)$. There exists a coupling
which reaches the minimum value. It is called a maximal coupling.
\end{Lemma}
\BLANC{ Corollary \ref{Cor1} is an immediate consequence of
Theorem \ref{Th} and Lemma \ref{lem_coupling}.
  Indeed, let $(Z_1,Z_2)$ be a maximal
coupling of $(\Dr_{\P_{\lambda}}({X_*}(t)),\mu)$. Combining
Theorem \ref{Th} and Lemma \ref{lem_coupling}, we obtain
$$
\P\left(Z_1\not=Z_2\right)\leq C e^{-\gamma t}\left(1+\int_{H}\abs{x}^2\,\lambda(dx)\right).
$$
Recall that, as explained in section $1.3$, $Z_2\in\H_{2}$ almost
surely. Hence
$$
\P\left(Z_1\not\in\H_{2}\right)\leq C e^{-\gamma t}\left(1+\int_{\H_2}\abs{x}^2\,\lambda(dx)\right).
$$
Since $\Dr(Z_1)=\Dr_{\P_{\lambda}}({X_*}(t))$, Corollary
\ref{Cor1} follows.}

 Let us first consider the case of the solutions of \eqref{galerkin}.
Assume that Hypothesis \ref{H0} holds.
 Let $N\in\N$ and $(x_0^1,x_0^2)\in\R^2$.
Combining arguments from \cite{KS1}, \cite{Matt}, it can be shown
 that there exists
 a decreasing function $p_N(\cdot)>0$ such that
\begin{equation}\label{toy3}
\norm{\left(\Pcal_1^N\right)^* \delta_{x_0^2}-\left(\Pcal_1^N\right)^* \delta_{x_0^1}}_{var}
\leq 1-p_N\left(\abs{x_0^1}+\abs{x_0^2}\right).
\end{equation}
Applying Lemma \ref{lem_coupling}, we build a maximal coupling $(Z_1,Z_2)=(Z_1(x_0^1,x_0^2),Z_2(x_0^1,x_0^2))$
of $
(\left(\Pcal_1^N\right)^* \delta_{x_0^1},\left(\Pcal_1^N\right)^* \delta_{x_0^2})$. It follows
\begin{equation}\label{toy4}
\P\left(Z_1=Z_2\right)\geq p_N\left(\abs{x_0^1}+\abs{x_0^2}\right)>0.
\end{equation}
Let $(W,\widetilde W)$ be a a couple of independent cylindrical
Wiener processes and $\delta>0$.
 We denote by $X_N(\cdot,x_0)$ and $\widetilde X_N(\cdot,x_0)$ the solutions
 of \eqref{galerkin} associated to $W$ and $\widetilde W$, respectively.
Now we build a couple of random variables $(V_1,V_2)=(V_1(x_0^1,x_0^2),V_2(x_0^1,x_0^2))$ on $P_NH$ as follows
\begin{equation}\label{toy5}
\Espace
(V_1,V_2)=
\left\{
\begin{array}{ll}
(X_N(\cdot,x_0),X_N(\cdot,x_0))& \textrm{ if } x_0^1=x_0^2=x_0,\\
(Z_1(x_0^1,x_0^2),Z_2(x_0^1,x_0^2))& \textrm{ if } (x_0^1,x_0^2)\in B_H(0,\delta)\backslash\{x_0^1=x_0^2\},\\
(X_N(\cdot,x_0^1),\widetilde X_N(\cdot,x_0^2))& \textrm{ else},
\end{array}
\right.
\end{equation}
where $B_H(0,\delta)$ is the ball of $H\times H$ with radius
$\delta$.

Then $(V_1(x_0^1,x_0^2),V_2(x_0^1,x_0^2))$ is a coupling of $
(\left(\Pcal_1^N\right)^* \delta_{x_0^1},\left(\Pcal_1^N\right)^*
\delta_{x_0^2})$. It can be shown that it depends measurably on
$(x_0^1,x_0^2)$. We then build a coupling $(X^1,X^2)$ of
$(\Dr(X_N(\cdot,x_0^1)),\Dr(X_N(\cdot,x_0^2)))$ by induction on
$\N$.
 We first set
$X^i(0)=x_0^i$ for $i=1,2$. Then, assuming that we have built
$(X^1,X^2)$ on $\{0,1,\dots,k\}$, we take $(V_1,V_2)$ as above
 independent of $(X^1,X^2)$ and set
$$
X^i(k+1)=V_i(X^1(k),X^2(k))\quad\textrm{ for } i=1,2.
$$
Taking into account \eqref{Eq1.10bis}, it is easily shown that the
time of return of $(X^1,X^2)$ in $B(0,4(c/\mu_1)B_1)$
 admits an exponential moment. We choose $\delta=4(c/\mu_1)B_1$. It follows from \eqref{toy4}, \eqref{toy5} that, $(X^1(n),X^2(n))\in B(0,\delta)$ implies that
the probability of having $(X^1,X^2)$ coupled (i.e. equal)  at
time $n+1$ is bounded below by $p_N(2\delta)>0$. Finally, remark
that if $(X^1,X^2)$ are coupled at time $n+1$, then they remain
coupled
 for any time after. Combining these three properties and using the fact that $(X^1(n),X^2(n))_{n\in\N}$ is a
 discrete strong Markov process, it is easily shown that
\begin{equation}\label{toy6}
\P\left(X^1(n)\not=X^2(n)\right)\leq C_Ne^{-\gamma_N
n}\left(1+\abs{x_0^1}^2+\abs{x_0^2}^2\right),
\end{equation}
with $\gamma_N>0$.

Recall that $(X^1,X^2)$ is a coupling of
$(\Dr(X_N(\cdot,x_0^1)),\Dr(X_N(\cdot,x_0^2)))$ on $\N$. It
follows
 that $(X^1(n),X^2(n))$ is a coupling of $((\Pcal_n^N)^* \delta_{x_0^1},(\Pcal_n^N)^* \delta_{x_0^2})$.
Combining Lemma \ref{lem_coupling} and \eqref{toy6}, we obtain, for $n\in\N$,
$$
\norm{\left(\Pcal_n^N\right)^* \delta_{x_0^2}-\left(\Pcal_n^N\right)^* \delta_{x_0^1}}_{var}\leq C_Ne^{-\gamma_N n}\left(1+\abs{x_0^1}^2+\abs{x_0^2}^2\right).
$$
Setting $n=\lfloor t\rfloor$
 and integrating $(x_0^2,x_0^1)$ over $((\Pcal_{t-n}^N)^*\lambda)\otimes\mu_N$ where $\mu_N$ is an invariant measure,
 it follows that, for any $\lambda\in P(P_NH)$,
\begin{equation}\label{toy7}
\norm{\left(\Pcal_t^N\right)^*\lambda-\mu_N}_{var}\leq C_N e^{-\gamma_N t}\left(1+\int_{P_N H}\abs{x}^2\,\lambda(dx)\right).
\end{equation}
This result is useless when considering equation \eqref{EqNS}
since the constants $C_N$, $\gamma_N$ strongly depend on $N$. If
one tries to apply directly the above arguments to the infinite
dimensional equation \eqref{EqNS}, one faces several difficulties.
First it is not known whether $\P_{x_0}$ is Markov. We only know
that, as explained in section $1.3$, a Markov transition
semi-group can be constructed. This is a major difficulty since
this property is implicitely used in numerous places above.
Another strong problem is that Girsanov transform is used in order
to obtain \eqref{toy3}. Contrary to the two dimensional case, no
Foias-Prodi estimate is available for the three dimensional
Navier-Stokes equations and the Girsanov transform should be done
in the infinite dimensional equation. This seems impossible. We
will show that we are able to prove an analogous result to
\eqref{toy3} by a completely different argument. However, this
will hold only for small initial data in $\H_2$. Another problem
will occur  since it is not known whether solutions starting in
$\H_2$ remain in $\H_2$.

We remedy the lack of Markov property by working only on Galerkin
approximations and prove that \eqref{toy7} holds with constants
uniform in $N$. As already mentioned, we prove that \eqref{toy3}
is true for $x_0^1$, $x_0^2$ in a small ball of $\H_2$ and
uniformly in $N$. Then, following the above argument, it remains
to prove that the time of return in this small ball admits an
exponential moment. Note that the smallness assumption on $f$ is
used at this step. In the following sections, we prove
\begin{Proposition}\label{ThN}
Assume that Hypothesis  \ref{H0} holds.  Then there exist
 ${{\delta^0}}={{\delta^0}}(B_0,D,\eps,\nu)$, $C=C(\phi,D,\eps,\nu)>0$ and $\gamma=\gamma(\phi,D,\eps,\nu)>0$
 such that if $\norm{f }_\eps^2\leq {{\delta^0}}\,$ holds, then, for any $N\in\N$, there exists
a unique invariant measure $\mu_N$ for
$\left(\Pcal_t^N\right)_{t\in \R^+}$. Moreover, for any $\lambda
\in \Pcal(P_N H)$
\begin{equation}\label{EqThN}
\norm{\left(\Pcal_t^N\right)^*\lambda-\mu_N}_{var}\leq C e^{-\gamma t}\left(1+\int_{P_N H}\abs{x}^2\,\lambda(dx)\right).
\end{equation}
\end{Proposition}
We now explain why this result implies Theorem \ref{Th}.

 Let $\lambda\in P(H)$ and $X_{\lambda}$ be a random variable on $H$ whose law is $\lambda$
 and which is independant of $W$.
 Since $\norm\cdot_{var}$ is the dual norm of $\abs\cdot_\infty$,
then  \eqref{EqThN} implies that
\begin{equation}\label{EqTh3}
\abs{\E\left(g(X_N(t,X_{\lambda}))\right)-\int_{P_NH}g(x)\,\mu_N(dx)}\leq
 C\abs{g}_\infty\left(1+\int_{H}\abs{x}^2\,\lambda(dx)\right),
\end{equation}
for any $g\in UC_b(\H_s)$ for $s<0$.

Assume that, for a subsequence $(N_k')_k$, $X_N(t,X_{\lambda})$
converges in distribution in $\H_s$
 to the law ${X_*}(t)$ under the weak solution $\P_\lambda$ of \eqref{EqNS}.
  Recall that the family $(\P_{\mu_N}^N)_N$ is tight.
 Hence, for a subsequence $(N_k)_k$ of $(N_k')_k$, $\P_{\mu_{N_k}}$ converges
to $\P_{\mu}$ a weak stationary solution of \eqref{EqNS} with initial law $\mu$.
 Taking the limit, \eqref{EqTh2} follows from \eqref{EqTh3},
 which yields Theorem \ref{Th}.

\section{Coupling of solutions starting from small initial data}

The aim of this section is to establish the following result. A
result analogous to \eqref{toy4} but uniform in $N$.
\begin{Proposition}\label{Prop_I_0}
Assume that Hypothesis \ref{H0} holds and that $f\in H$. Then
there exist $(T,\delta)\in(0,1)^2$ such that, for any $N\in\N$,
there exists a coupling $(Z_1(x_0^1,x_0^2),Z_2(x_0^1,x_0^2))$ of
$((\Pcal_T^N)^*\delta_{x_0^1},(\Pcal_T^N)^*\delta_{x_0^2})$ which
measurably depends on $(x_0^1,x_0^2)\in\H_2$ and which verifies
 \begin{equation}\label{I_0bis}
\P\left(Z_1(x_0^1,x_0^2)=Z_2(x_0^1,x_0^2)\right)\geq\frac{3}{4}
\end{equation}
 provided
\begin{equation}\label{2.5}
\norm{x_0^1}_2^2\vee\norm{x_0^2}_2^2\leq \delta.
 \end{equation}
\end{Proposition}
Assume that Hypothesis \ref{H0} holds  and that $f\in H$. Let
$T\in(0,1)$. Applying Lemma \ref{lem_coupling}, we build
$(Z_1(x_0^1,x_0^2),Z_2(x_0^1,x_0^2))$
 as the maximal coupling of $(\Pcal_T^*\delta_{x_0^1},\Pcal_T^*\delta_{x_0^2})$.
Measurable dependance follows from a slight extension of Lemma
$1.17$ (see \cite{ODASSO1}, remark $A.1$).

In order to establish Proposition \ref{Prop_I_0}, it is sufficient
to prove that
 there exists $c(B_0,D)$ not depending on $T\in(0,1)$ and on $N\in\N$
 such that
\begin{equation}\label{I_0ter}
 \norm{\left(\Pcal_T^N\right)^*\delta_{x_0^2}-\left(\Pcal_T^N\right)^*\delta_{x_0^1}}_{var}\leq c(B_0,D)\sqrt
 T,
\end{equation}
 provided
\begin{equation}\label{2.5ter}
\norm{x_0^1}_2^2\vee\norm{x_0^2}_2^2\leq B_0 T^3.
\end{equation}
Then it suffices to choose $T\leq 1/(4c(B_0,D))^2$ and
$\delta=B_0T^3$.

Since $\norm\cdot_{var}$ is the dual norm of $\abs\cdot_\infty$,
 \eqref{I_0ter} is equivalent to
\begin{equation}\label{I_0}
 \abs{\E\left( g(X_N(T,x_0^2))-g(X_N(T,x_0^1))\right)}\leq
8\abs{g}_\infty c(B_0,D)\sqrt T.
\end{equation}
 for any $g\in UC_b(P_N H)$.

 It follows from the density of $C^1_b(P_N H)\subset UC_b(P_N H)$
 that, in order
 to establish Proposition \ref{Prop_I_0}, it is sufficient to prove that \eqref{I_0} holds
for any $N\in\N$, $T\in(0,1)$ and $g\in C_b^1(P_N H)$ provided \eqref{2.5ter} holds.

The proof of \eqref{I_0} under this condition is splitted into the
next three subsections.

\subsection{A priori estimate}

\

\noindent For any process $X$, we define the $\H_1$--energy of $X$ at time $t$ by
$$
 E_X^{\H_1}(t)=\norm{X(t)}^2+ \int_{0}^t\norm{X(s)}_2^2ds.
$$
Now we establish the following result which will be useful in the proof of \ref{I_0}.
\begin{Lemma}\label{Prop_nrj}
Assume that Hypothesis \ref{H0} holds  and that $f\in H$. There
exist $K_0=K_0(D)$ and $c=c(D)$
 such that for any $T\leq 1$ and any  $N\in\N$,
 we have
$$
\P\left(\sup_{(0,T)} E_{X_N(\cdot,x_0)}^{\H_1} >K_0\right)\leq c\left(1+\frac{B_0}{K_0}\right)\sqrt T,
$$
provided $\norm{x_0}^2\leq B_0 T$.
\end{Lemma}

\noindent Let $X_N=X_N(\cdot,x_0)$. Ito Formula on
$\norm{{X_N}}^2$ gives
\begin{equation}\label{8.1}
d\norm{{X_N}}^2+2  \norm{{X_N}}_2^2dt=dM_{\H_1}+I_{\H_1}dt+\norm{{P_N \phi}({X_N})}^2_{\mathcal L_2(H;\H_1)}dt
+I_fdt,
\end{equation}
where
$$
\Espace
\left\{
\begin{array}{l}
I_{\H_1}=-2\left(A{X_N},B({X_N})\right),\;I_f=2\left(A{X_N},f\right),\\
\;M_{\H_1}(t)=2\int_0^t\left(A{X_N}(s),\phi({X_N}(s))dW(s)\right).
\end{array}
\right.
$$
Combining a H\"older inequality, a Agmon inequality and a arithmetico-geometric inequality gives
\begin{equation}\label{8.2}
I_{\H_1}\leq 2\norm{{X_N}}_2\abs{{X_N}}_\infty\norm{{X_N}}\leq
c\norm{{X_N}}_2^\frac{3}{2}\norm{{X_N}}^\frac{3}{2} \leq
\frac{1}{4}\norm{{X_N}}_2^2+c\norm{{X_N}}^6.
\end{equation}
Similarly, using Poincar\'e inequality and Hypothesis \ref{H0},
\begin{equation}\label{8.2bis}
I_{f}\leq \frac{1}{4}\norm{{X_N}}_2^2+c\abs{f}^2\leq
\frac{1}{4}\norm{{X_N}}_2^2+cB_0.
\end{equation}
We deduce from \eqref{8.1}, \eqref{8.2}, \eqref{8.2bis},
Hypothesis \ref{H0} and Poincar\'e inequality that
\begin{equation}\label{8.3}
d\norm{{X_N}}^2+  \norm{{X_N}}_2^2dt\leq dM_{\H_1}+cB_0 dt
+c\norm{{X_N}}^2 \left(\norm{{X_N}}^4-4K_0^2\right)\,dt,
\end{equation}
where
\begin{equation}\label{2.9bis}
K_0=\sqrt{\frac{\mu_1}{8c}}.
\end{equation}
 \noindent Setting
$$
\sig_{\H_1}=\inf\left\{t\in(0,T)\,\left|\,  \norm{{X_N}(t)}^2> 2 K_0   \right.\right\},
$$
we infer from $\norm{x_0}^2\leq B_0 T$ that for any $t\in(0,\sig_{\H_1})$
\begin{equation}\label{8.4}
E_{X_N}^{\H_1}(t)\leq cB_0T +M_{\H_1}(t).
\end{equation}
We deduce from  Hypothesis \ref{H0} and from Poincar\'e inequality that $\phi(x)^*A$ is bounded
in $\mathcal L(\H_1;\H_1)$ by $cB_0$. It follows that for any $t\in (0,\sig_{\H_1})$
$$
\left<M_{\H_1}\right>(t)=4\int_0^t\norm{{P_N
\phi}({X_N}(s))^*A{X_N}(s)}^2dt\leq
 cB_0 \int_0^t\norm{{X_N}(s)}^2ds\leq 2cK_0B_0 T.
$$
Hence a Burkholder-Davis-Gundy inequality gives
$$
\E\left(\sup_{(0,\sig_{\H_1})}M_{\H_1}\right)\leq
c\E\sqrt{\left<M_{\H_1}\right>(\sig_{\H_1})}\leq
c\sqrt{K_0B_0T}\leq c(K_0+B_0)\sqrt T.
$$
It follows from \eqref{8.4} and $T\leq 1$ that
$$
\E\left(\sup_{(0,\sig_{\H_1})}E_{X_N}^{\H_1}\right)\leq
c(B_0+K_0)\sqrt T,
$$
which yields, by a Chebyshev inequality,
$$
\P\left(\sup_{(0,\sig_{\H_1})}E_{X_N}^{\H_1}> K_0\right)\leq c\left(1+\frac{B_0}{K_0}\right)\sqrt T.
$$
Now, since $\sup_{(0,\sig_{\H_1})}E_{X_N}^{\H_1}\leq K_0$ implies $\sig_{\H_1}=T$, we deduce
Lemma \ref{Prop_nrj}.

\subsection{Estimate of the derivative of $X_N$}\label{Preuve_eta}
\

 Let $N\in\N$ and $(x_0,h)\in (\H_2)^2$. 
We are concerned with the following equation
\begin{equation}\label{Eqeta}
\Espace
\left\{
\begin{array}{rcl}
d{\eta_N}+   A {\eta_N} \;dt+{P_N \widetilde B}({X_N},{\eta_N})\, dt &=& P_N({ \phi}'({X_N})\cdot{\eta_N})\, dW,\\
                             {\eta_N}(s,s,x_0)\cdot h &=& P_N h,
\end{array}
\right.
\end{equation}
where $\widetilde
B({X_N},{\eta_N})=B({X_N},{\eta_N})+B({\eta_N},{X_N})$,
$X_N=X_N(\cdot,x_0)$ and ${\eta_N}(t)={\eta_N}(t,s,x_0)\cdot h$
for $t\geq s$.

 Existence and uniqueness of the solutions of \eqref{Eqeta} are easily shown.
 Moreover if $g\in C_b^1(P_N H)$, then, for any $t\geq 0$, we have
\begin{equation}\label{Eqetabis}
\left(\nabla\left(\Pcal_t^Ng\right)(x_0),h\right)=\E\left(\nabla g(X_N(t,x_0)),\eta_N(t,0,x_0)\cdot h\right).
\end{equation}
For any process $X$, we set
\begin{equation}\label{2.12bis}
\sig(X)=\inf\left\{t\in
(0,T)\,\left|\,\int_0^t\norm{X(s)}_2^2ds\geq K_0+1\right.\right\},
\end{equation}
where $K_0$ is defined in Lemma \ref{Prop_nrj}. We establish the
following result.
\begin{Lemma}\label{Prop_eta}
Assume that  Hypothesis \ref{H0}  holds  and that $f\in H$. Then
there exists $c=c(B_0,D)$ such that for any $N\in \N$, $T\leq 1$
and $(x_0,h)\in (\H_2)^2$
$$
\E\int_0^{\sig(X_N(\cdot,x_0))} \norm{{\eta_N}(t,0,x_0)\cdot h}_3^2 \,dt\leq
 c\norm{h}_2^2.
$$
\end{Lemma}

\noindent For a better readability, we set ${\eta_N}(t)={\eta_N}(t,0,x_0)\cdot h$ and $\sig=\sig(X_N(\cdot,x_0))$.
 Ito Formula on $\norm{{\eta_N}(t)}_2^2$
gives
\begin{equation}\label{5.1}
d\norm{{\eta_N}}_2^2+2 \norm{{\eta_N}}_3^2dt=dM_{{\eta_N}}+I_{\eta_N} \,dt+
\norm{P_N \left(\phi'({X_N})\cdot {\eta_N}\right)}^2_{\mathcal L_2(U;\H_2)}dt,
\end{equation}
where
$$
\Espace
\left\{
\begin{array}{rcl}
M_{{\eta_N}}(t)&=&2\int_0^t \left(A^2{\eta_N},({P_N \phi}'({X_N})\cdot {\eta_N})\, dW\right)\,ds,\\
I_{\eta_N}&=&-2\left(A^\frac{3}{2}{\eta_N},A^\frac{1}{2}\widetilde
B({X_N},{\eta_N})\right).
\end{array}
\right.
$$
It follows from H\"older inequalities, Sobolev Embedding and a
arithmetico-geometric inequality
$$
I_{\eta_N}\leq c \norm{{\eta_N}}_3\norm{{\eta_N}}_2\norm{{X_N}}_2\leq  \norm{{\eta_N}}_3^2+c\norm{{\eta_N}}_2^2\norm{{X_N}}_2^2.
$$
Hence, we deduce from \eqref{5.1} and   Hypothesis \ref{H0}
$$
d\norm{{\eta_N}}_2^2+ \norm{{\eta_N}}_3^2dt\leq dM_{{\eta_N}}+c\norm{{\eta_N}}_2^2\norm{{X_N}}_2^2+B_0\norm{{\eta_N}}^2_2dt.
$$
Integrating and taking the expectation, we obtain
\begin{equation}\label{5.2}
\E\left( \mathcal E(\sig,0)\norm{{\eta_N}(\sig)}_2^2
+ \int_0^\sig \mathcal E(\sig,t)\norm{{\eta_N}(t)}_3^2\,dt\right)
\leq\norm{h}^2_2,
\end{equation}
where
$$
\mathcal E(t,s)=e^{-B_0t-c\int_s^t\norm{{X_N}(r)}_2^2\,dr}.
$$

\noindent Applying the definition of $\sig$, we deduce
\begin{equation}\label{5.3}
\E\int_0^\sig \norm{{\eta_N}(t)}_3^2 \,dt\leq
 \norm{h}_2^2\exp\left(c(K_0+1)+B_0 T\right),
\end{equation}
which yields Lemma \ref{Prop_eta}.

\subsection{Proof of \eqref{I_0}}
\

 \noindent Let $\psi\in C^\infty\left(\R; [0,1]\right)$ such that
$$
\psi=0 \;\textrm{ on } \;(K_0 +1,\infty)\quad \textrm{ and }\quad \psi=1 \;\textrm{ on }\; (-\infty,K_0 ).
$$
For any process $X$, we set
$$
 \psi_{X}=\psi\left(\int_0^{T}\norm{{X}(s)}_2^2ds\right).
$$
\noindent Remark that
\begin{equation}\label{2.10}
\abs{\E\left( g({X_N}(T,x_0^2)) -g({X_N}(T,x_0^1))\right)}\leq I_0+\abs{g}_\infty\left(I_1+I_2\right),
\end{equation}
where
$$
\Espace
\left\{
\begin{array}{lcl}
I_0 &=& \abs{\E\left(g({X_N}(T,x_0^2))\psi_{{X_N(\cdot,x_0^2)}}-g({X_N}(T,x_0^1))\psi_{{X_N(\cdot,x_0^1)}}\right)},\\
I_i &=&\P\left(\int_0^T\norm{{X_N}(s,x_0^i)}_2^2ds>K_0 \right).
\end{array}
\right.
$$
For any $\theta\in[1,2]$, we set
$$
\Espace
\left\{
\begin{array}{ll}
x_0^\theta=(2-\theta)x_0^1+(\theta-1)x_0^2,& X_\theta={X_N}(\cdot,x_0^\theta),\\
\eta_\theta(t)=\eta_N(t,0,x_0^\theta),&\sig_\theta=\sig(X_\theta).
\end{array}
\right.
$$
Recall that $\sig$ was defined in \eqref{2.12bis}. For a better
readability, the dependance on $N$ has been omitted.
 \noindent Setting
$$
h=x_0^2-x_0^1,
$$
we have
\begin{equation}\label{2.11}
I_0\leq \int_1^2\abs{J_\theta}d\theta\quad
J_\theta=\left(\nabla \E\left(g(X_{\theta}(T))\psi_{X_\theta}\right),h\right).
\end{equation}
 To bound $J_\theta$, we apply a truncated Bismut-Elworthy-Li formula (See appendix \ref{Preuve_2.15})
\begin{equation}\label{2.15}
J_\theta=\frac{1}{T}J_{\theta,1}'+2J_{\theta,2}',
\end{equation}
where
$$
\Espace
\left\{
\begin{array}{lcl}
J_{\theta,1}'&=&\E \left(g(X_\theta(T))\psi_{X_\theta}\int_0^{\sig_\theta} (\phi^{-1}(X_\theta(t))\cdot \eta_\theta(t)\cdot h,dW(t))\right),\\
J_{\theta,2}'&=&\E\left(g(X_\theta(T))\psi'_{X_\theta} \int_0^{\sig_\theta} \left(1-\frac{t}{T}\right)\left(AX_\theta(t),A (\eta_\theta(t)\cdot h)\right)\,dt\right),\\
\psi'_X&=&\psi'\left(\int_0^{T}\norm{X_\theta(s)}_2^2ds\right).
\end{array}
\right.
$$
 It follows from H\"older inequality that
$$
\abs{J_{\theta,2}'}\leq \abs{g}_\infty\abs{\psi'}_\infty \sqrt{\E\int_0^{\sig_\theta} \norm{X_\theta(t)}_2^2 \,dt }
\sqrt{\E\int_0^{\sig_\theta} \norm{\eta_\theta(t)\cdot h}_2^2 \,dt }.
$$
and from Hypothesis \ref{H0} that
$$
\abs{J_{\theta,1}'}\leq \abs{g}_\infty B_0\sqrt{\E\int_0^{\sig_\theta} \norm{\eta_\theta(t)\cdot h}_3^2 \,dt }.
$$
Hence for any $T\leq 1$
\begin{equation}\label{2.16}
\abs{J_\theta}\leq c(B_0,D)\abs{g}_\infty \frac{1}{T}\sqrt{\E\int_0^{\sig_\theta} \norm{\eta_\theta(t)\cdot h}_3^2 \,dt }.
\end{equation}
\noindent Combining \eqref{2.16} and Lemma \ref{Prop_eta}, we obtain
$$
\abs{J_\theta}\leq c(B_0,D)\abs g_\infty\frac{\norm{h}_2}{T},
$$
which yields, by \eqref{2.5ter} and \eqref{2.11},
$$
I_0\leq c(B_0,D)\abs g_\infty\sqrt{T}.
$$
\noindent Since $B_0T^3\leq B_0T$, we can apply Lemma
\ref{Prop_nrj} to control $I_1+I_2$ in \eqref{2.10} if
\eqref{2.5ter} holds. Hence \eqref{I_0} follows provided
\eqref{2.5ter} holds, which yields Proposition \ref{Prop_I_0}.

\section{Time of return in a small ball of $\H_2$}\label{Preuve_Lyap}

Assume that Hypothesis \ref{H0} holds. Let $N\in\N$ and  $T,\delta,Z_1,Z_2$ be as in Proposition \ref{Prop_I_0}.
Let $(W,\widetilde W)$ be a couple of independant cylindrical Wiener processes on $H$.
 We denote by $X_N(\cdot,x_0)$ and $\widetilde X_N(\cdot,x_0)$ the solutions
 of \eqref{galerkin} associated to $W$ and $\widetilde W$, respectively.
We build a couple of random variables $(V_1,V_2)=(V_1(x_0^1,x_0^2),V_2(x_0^1,x_0^2))$ on $P_NH$ as follows
\begin{equation}\label{toy5bisbis}
\Espace
(V_1,V_2)=
\left\{
\begin{array}{ll}
(X_N(\cdot,x_0),X_N(\cdot,x_0))& \textrm{ if } x_0^1=x_0^2=x_0,\\
(Z_1(x_0^1,x_0^2),Z_2(x_0^1,x_0^2))& \textrm{ if } (x_0^1,x_0^2)\in B_{\H_2}(0,\delta)\backslash\{x_0^1=x_0^2\},\\
(X_N(\cdot,x_0^1),\widetilde X_N(\cdot,x_0^2))& \textrm{ else},
\end{array}
\right.
\end{equation}
 We then build $(X^1,X^2)$ by induction on $T\N$.
 Indeed, we first set
$X^i(0)=x_0^i$ for $i=1,2$. Then, assuming that we have built
$(X^1,X^2)$ on $\{0,T,2T,\dots,nT\}$,
 we take $(V_1,V_2)$ as above independent of $(X^1,X^2)$ and we set
$$
X^i((n+1)T)=V_i(X^1(nT),X^2(nT))\quad\textrm{ for } i=1,2.
$$
It follows that $(X^1,X^2)$ is a discrete strong Markov process
and a coupling of $(\Dr(X_N(\cdot,x_0^1)),\Dr(X_N(\cdot,x_0^2)))$
on $T\N$. Moreover, if $(X^1,X^2)$ are coupled at time $nT$, then
they remain coupled for any time after.

We set
\begin{equation}\label{3.1bis}
\tau=\inf\left\{t\in
T\N\backslash\{0\}\,\left|\,
\norm{X^1(t)}_2^2\vee\norm{X^2(t)}_2^2\leq \delta \right.\right\}.
\end{equation}
The aim of this section is to establish the following result.
\begin{Proposition}\label{Prop_Lyap}
Assume that Hypothesis \ref{H0} holds. There exist
 $ \delta^3= \delta^3(B_0,D,\eps,\delta)$, $\alpha=\alpha(\phi,D,\eps,\delta)>0$ and
$K"=K"(\phi,D,\eps,\delta)$ such that for any $(x_0^1,x_0^2)\in H\times H$ and any $N\in\N$
$$
\E\left(e^{\alpha \tau}\right)\leq K"\left(1+\abs{x_0^1}^2+\abs{x_0^2}^2\right),
$$
 provided $\norm f^2_\eps \leq \delta^3$.
\end{Proposition}

The result is based on the fact that, in the absence of noise and
forcing term, all solutions go to zero exponentially fast in $H$.
A similar idea is used for the two-dimensional Navier-Stokes
equations in \cite{KS1}. The proof is based on the following four
Lemmas. The first one allows to control the probability that the
contribution of the noise is small. Its proof strongly uses the
assumption that the noise is diagonal in the eigenbasis of $A$. As
already mentioned, in the additive case, the proof is easy and
does not need this assumption.
\begin{Lemma}\label{PropSmallH3}
Assume that  Hypothesis \ref{H0} holds. For any $t, {M} >0$, there
exists
 $p_0(t, {M} )=p_0(t, {M} ,\eps,(\abs{\phi_n}_\infty)_n,D)>0$
 such that for any adapted process $X$
$$
\P\left(\sup_{(0,t)}\norm{Z}_2^2\leq  {M}  \right)\geq p_0(t, {M}
),
$$
where
$$
Z(t)=\int_0^te^{-A(t-s)}\phi(X(s))dW(s).
$$
\end{Lemma}
It is proved in section $3.1$.

Then, using this estimate and the smallness assumption on the
forcing term, we estimate the moment of the first return time in a
small ball in $H$.




\noindent Let $\delta_3>0$. We set
$$
\tau_{L^2}=\tau\wedge\inf\left\{t\in T\N^*\,\left|\,   \abs{X^1(t)}^2\vee\abs{X^2(t)}\geq \delta_3   \right.\right\}.
$$
 \begin{Lemma}\label{lem7.2}
Assume that  Hypothesis \ref{H0} holds. Then, for any
$\delta_3>0$, there exist $C_3(\delta_3)$, $C_3'(\delta_3)$ and
 $\gamma_3(\delta_3)$
 such that for any $(x_0^1,x_0^2)\in (\H_2)^2$
$$
 \E\left(e^{\gamma_3
\tau_{L^2}}\right)\leq
C_3\left(1+\abs{x_0^1}^2+\abs{x_0^2}^2\right), $$
provided
$$
\abs{f}\leq C_3'.
$$
\end{Lemma}
 \noindent The proof is postponed to section
\ref{Preuve7.2}.

Then, we need to get a finer estimate in order to control the time
necessary to enter a ball in stronger topologies. To prove the two
next lemmas, we use an argument similar to one used in the
determinist theory (see \cite{Temam2}, chapter $7$).

\begin{Lemma}\label{lem7.4}
Assume that  Hypothesis \ref{H0} holds. Then, for any $\delta_4$,
there exist $p_4(\delta_4)>0$, $C_4'(\delta_4)>0$ and
$R_4(\delta_4)>0$
 such that for any $x_0$ verifying $\abs{x_0}^2\leq R_4$, we have for any $T\leq 1$
$$
\P\left(\norm{X_N(T,x_0)}^2\leq \delta_4\right)\geq p_4,
$$
provided
$$
\abs{f}\leq C_4'.
$$
\end{Lemma}

\noindent The proof  is postponed
 to section \ref{Preuve7.4second}.

\begin{Lemma}\label{lem7.5}
Assume that  Hypothesis \ref{H0} holds. Then, for any $\delta_5$,
there exist $p_5(\delta_5)>0$, $C_5'(\delta_5)>0$ and
$R_5(\delta_5)>0$
 such that for any $x_0$ verifying $\norm{x_0}^2\leq R_5$   and  for any $T\leq 1$
$$
 \P\left(\norm{X_N(T,x_0)}_2^2\leq
\delta_5\right)\geq p_5.
$$
provided
$$
\norm{f}_\eps\leq C_5'.
$$
\end{Lemma}

\noindent The proof  is postponed to section \ref{Preuve7.5}.

\

{\it Proof of Proposition \ref{Prop_Lyap}}: We set
$$
\delta_5=\delta,\quad \delta_4=R_5(\delta_5),\quad
\delta_3=R_4(\delta_4),\quad p_4=p_4(\delta_4), \quad
p_5=p_5(\delta_5), \quad p_1=\left(p_4p_5\right)^ 2,
$$
and
$$
\delta^3=C_3'(\delta_3)\wedge C_4'(\delta_4)\wedge C_5'(\delta_5).
$$
 By the definition of $\tau_{L^2}$, we have
$$
\abs{X^1(\tau_{L^2})}^2\vee\abs{X^2(\tau_{L^2})}^2\leq
R_4(\delta_4).
$$
 We distinguish
three cases.

The first case is $\norm{X^1(\tau_{L^2})}_2^2\vee\norm{X^2(\tau_{L^2})}_2^2\leq\delta$, which obviously yields
\begin{equation}\label{7.5bis1}
\P\left(\min_{k=0,\dots,2}\max_{i=1,2}\norm{X^i(\tau_{L^2}+kT)}_2^2\leq\delta\,\left|\,
\left(X^2(\tau_{L^2}),X^2(\tau_{L^2})\right)\right.\right)\geq
p_1.
\end{equation}

We now treat the case $x_0=X^1(\tau_{L^2})=X^2(\tau_{L_2})$ with
$\norm{x_0}_ 2^2>\delta$. Combining  Lemma \ref{lem7.4} and Lemma
\ref{lem7.5}, we deduce from the weak Markov property of $X_N$
that
$$
\P\left(\norm{X_N(2T,x_0)}_2^2\leq \delta\right)\geq p_5p_4,
$$
provided $\abs{x_0}^2\leq R_4$.
Recall  that, in that case, $X^1(\tau_{L^2}+2T)=X^2(\tau_{L^2}+2T)$.
 Hence, since the law of $X^1(\tau_{L^2}+2T)$ conditioned by $(X^1(\tau_{L^2}),X^2(\tau_{L_2}))$ is
$\Dr(X_N(2T,x_0))$, it follows
$$
\P\left(\max_{i=1,2}\norm{X^i(\tau_{L^2}+2T)}_2^2\leq\delta\,\left|\,
\left(X^1(\tau_{L^2}),X^2(\tau_{L_2})\right)\right.\right)\geq
p_4p_5\geq p_1,
$$
and then \eqref{7.5bis1}

The last case is $X^1(\tau_{L^2})\not=X^2(\tau_{L_2})$ and
$\norm{X^1(\tau_{L^2})}_2^2\vee\norm{X^2(\tau_{L^2})}_2^2>\delta$. In that case, $(X^1(\tau_{L^2}+T),X^2(\tau_{L^2}+T))$
 conditioned by $(X^1(\tau_{L^2}),X^2(\tau_{L_2}))$ are independent. Hence,
since the law of $X^i(\tau_{L^2}+T)$ conditioned by
$(X^1(\tau_{L^2}),X^2(\tau_{L_2}))=(x_0^1,x_0^2)$ is
 $\Dr(X_N(T,x_0^i))$, it follows from Lemma \ref{lem7.4} that
$$
\P\left(\max_{i=1,2}\norm{X^i(\tau_{L^2}+T)}_1^2\leq\delta_4\,\left|\,
\left(X^1(\tau_{L^2}),X^2(\tau_{L_2})\right)\right.\right)\geq
p_4^2.
$$
Then, we distinguish the three cases
$(\norm{X^i(\tau_{L^2}+T)}_1^2)_{i=1,2}$ in the small ball of
$\H_2$, equal or different and we deduce from Lemma \ref{lem7.5}
by the same method
$$
\P\left(\min_{k=1,2}\max_{i=1,2}\norm{X^i(\tau_{L^2}+kT)}_2^2\leq\delta\,\left|\,
\left(X^1(\tau_{L^2}+T),X^2(\tau_{L_2}+T)\right)\right.\right)\geq
p_5^2,
$$
provided
$$
\max_{i=1,2}\norm{X^i(\tau_{L^2}+T)}_1^2\leq\delta_4.
$$
Combining the two previous inequalities, we deduce \eqref{7.5bis1}
for the latter case. We have thus proved that \eqref{7.5bis1} is
true almost surely.

Integrating \eqref{7.5bis1}, we obtain
\begin{equation}\label{7.5bis}
\P\left(\min_{k=0,\dots,
2}\max_{i=1,2}\norm{X^i(\tau_{L^2}+kT)}_2^2\leq\delta\right)\geq
p_1.
\end{equation}

\noindent Combining Lemma \ref{lem7.2} and \eqref{7.5bis}, we
conclude.

\subsection{Probability of having a small noise}
\


We now establish Lemma \ref{PropSmallH3}.

\noindent We deduce from H\"older inequality and from
$\sum_n\mu_n^{-2}<\infty$ that Hypothesis \ref{H0} implies the
following fact: for any $\eps_0\in (0,\eps)$, there exists $\alpha
\in (0,1)$, a family $(\bar \phi_n)_n$ of measurable maps $H\to
\R$  and a family $(b_i)_i$ of positive numbers  such that
\begin{equation}\label{3.6bis}
\Espace \left\{
\begin{array}{l}
 \phi(x)\cdot e_n =b_n\bar \phi_n(x)e_n,\\
 \sup_{x\in H}\abs{\bar \phi_n(x)}\leq 1,
\quad B^*=\sum_n\,\mu_n^{1+\eps_0}\,\left( b_n\right)^{2(1-\alpha)}<\infty.
\end{array}
\right.
\end{equation}
\noindent For simplicity we restrict our attention to the case
$t=1$. The generalization is easy.

\noindent Remark that
$$
Z(t)=\sum_n b_n Z_n(t)e_n,
$$
where
$$
 Z_n(t)=\int_0^te^{- \mu_n(t-s)}\bar \phi_n(X(s))dW_n(s),\;\textrm{ where }\;W=\sum_nW_ne_n.
$$
It follows from $\norm{Z}^2_2=\sum_n b_n^2 \mu_n \abs{\sqrt{
\mu_n}Z_n}^2$ and from \eqref{3.6bis} that
\begin{equation}\label{Small3.1}
\P\left(\sup_{(0,1)}\norm{Z}^2_2 \leq B^* {M} \right)\geq
\P\left(\sup_{(0,1)}\abs{\sqrt{ \mu_n}Z_n}^2 \leq  {M}  \,\mu_n^{\eps_0}\left( b_n\right)^{-2\alpha} ,
\;\forall\; n\right).
\end{equation}
Setting
$$
W_n'(t)=\sqrt{ \mu_n}W_n\left(\frac{t}{ \mu_n}\right),
$$
we obtain $(W_n')_n$ a family   of independent brownian motions.
Moreover we have
$$
\sqrt{ \mu_n}Z_n(t)=Z_n'( \mu_n t),
$$
where
$$
Z_n'(t)=\int_0^te^{-(t-s)}\psi_n(s)dW_n'(s), \quad \psi_n(s)=\bar \phi_n\left(X\left(\frac{s}{  \mu_n}\right)\right).
$$
Hence, it follows from \eqref{Small3.1} that
\begin{equation}\label{Small3.3}
\P\left( \sup_{(0,1)}\norm{Z}^2_2 \leq B^* {M} \right)\geq
\P\left(\sup_{(0, \mu_n)}\abs{Z_n'}^2 \leq  {M}  \,\mu_n^{\eps_0}\left( b_n\right)^{-2\alpha} ,
\;\forall\; n\right).
\end{equation}
Let $W_{n,i}'=W_n'(i+\cdot)-W_n'(i)$ on $(0,1)$. We set
$$
\Espace
M_{n,i}(t)=
\left\{
\begin{array}{lcl}
0 &\textrm{ if }& t\leq 0,\\
\int_0^{1\wedge t}e^{s}\,\psi_n(i+s)dW_{n,i}'(s)&\textrm{ if }& t\geq 0.
\end{array}
\right.
$$
Remark that
$$
Z_n'(t)=\sum_{i=1}^\infty e^{-(t-i)}M_{n,i}(t-i),
$$
which yields for any $q\in\N$
\begin{equation}\label{Small3.4}
\sup_{(0,q)}\norm{Z_n'}_2\leq \left(\frac{e}{e-1}\right)\max_{i=0,...,q-1}\sup_{(0,1)}\abs{M_{n,i}}.
\end{equation}
Remark that $(W_{n,k}')_{n,k}$ is a family of independant brownian motions on $(0,1)$. It follows that
$(M_{n,k})_{n,k}$ are martingales verifying $\left<M_{n,k},M_{n',k'}\right>=0$ if $(n,k)\not =(n',k')$.

\noindent Hence, combining a Theorem by Dambis, Dubins and
Schwartz (Theorem 4.6 page 174 of \cite{KARATZAS}) and a Theorem
by Knight (Theorem 4.13 page 179 of \cite{KARATZAS}), we obtain a
family $(B_{n,k})_{n,k}$  of independent brownian motions
verifying
\begin{equation}\label{Small3.5}
M_{n,k}(t)=B_{n,k}(\left<M_{n,k}\right>(t)).
\end{equation}
\begin{Remark}
In the two previous Theorem (Theorem 4.6 page 174 and Theorem 4.13
page 179 of \cite{KARATZAS}), it is assumed that
$\P\left(<M>(\infty)<\infty\right)=0.$ However, as explained in
Problem 4.7 of \cite{KARATZAS}, the proof is easily adapted  to
the case $\P\left(<M>(\infty)<\infty\right)>0.$
\end{Remark}
Remarking that for any $t\in (0,1)$
$$
\left<M_{n,k}\right>(t)=\int_0^{t}\abs{\psi_n(k+s)}^2ds\leq 1,
$$
 we deduce from \eqref{Small3.4} and \eqref{Small3.5} that for any $q\in\N^*$
$$
\sup_{(0,q)}\norm{Z_n'}_2\leq \left(\frac{e}{e-1}\right)\max_{i=0,...,q-1}\sup_{(0,1)}\abs{B_{n,i}}.
$$
Hence it follows from \eqref{Small3.3} that
$$
\P\left(\sup_{(0,1)}\norm{Z}^2_2 \leq  {M} \right)\geq
\P\left(\sup_{(0,1)}\abs{B_{n,i}}^2\leq c {M} \,\mu_n^{\eps_0}\left( b_n\right)^{-2\alpha},\;
\forall \,n,\;\forall \,i\leq \mu_n+1\right),
$$
where $c=\left(\frac{e-1}{  e}\right)^2\frac{1}{B^*}$.

\noindent We deduce from the independence of $(B_{n,k})_{n,k}$
that
\begin{equation}\label{Small3.6}
\P\left(\sup_{(0,1)}\norm{Z}^2_2 \leq  {M} \right)\geq
\prod_{n\in\N^*}\left( P\left( c  {M} \,\mu_n^{\eps_0}\left( b_n\right)^{-2\alpha}\right)^{ \mu_n+1}\right),
\end{equation}
where
$$
P(d_0)=\P\left(\sup_{(0,1)}\abs{B_{1,1}}^2 \leq d_0 \right).
$$

\noindent Recall there exists a family $(c_p)_p$ such that
$$
\E\left(\sup_{(0,1)}\abs{B_{1,1}}^{2p}\right)\leq c_p.
$$
It follows from Chebyshev inequality and from $1-x\geq e^{-ex}$ for any $x\leq e^{-1}$
 that for any $d_0\leq d_p=\left(e^{-1}c_p\right)^\frac{1}{p}$
$$
P(d_0)\geq 1-c_pd_0^{-p}\geq e^{-ec_pd_0^{-p}}.
$$
Applying \eqref{Small3.6}, we obtain for any $p>0$
\begin{equation}\label{Small3.7}
\P\left(\sup_{(0,1)}\norm{Z}^2_2 \leq  {M} \right)\geq
C_p( {M} )\exp\left(-\frac{c_p'}{ {M} ^p}
\sum_{n> \N(p, {M} )}\left(\frac{ \mu_n+1}{\mu_n^{\eps_0p}}\right)b_n^{2\alpha p}
\right),
\end{equation}
where
$$
\Espace
\left\{
\begin{array}{lcl}
N(p, {M} )&=&\sup\left\{n\in\N\backslash\{0\}\,\left|\,  {M} \,\mu_n^{\eps_0}\left( b_n\right)^{-2\alpha}
 \leq d_p \right.\right\},\\
C_p( {M} )
&=&\Pi_{n\leq\N(p, {M} )}\left( P\left(c {M} \,\mu_n^{\eps_0}\left( b_n\right)^{-2\alpha}\right)^{ \mu_n+1}\right).
\end{array}
\right.
$$
Choosing $p$ sufficiently high, we deduce from {\bf H0}
that
$$
\sum_{n}\left(\frac{ \mu_n+1}{\mu_n^{\eps_0p}}\right)b_n^{2\alpha p}\leq C_p'<\infty,
$$
which yields, by \eqref{Small3.7}, that for any $ {M} >0$ and for
$p$ sufficiently high
\begin{equation}\label{Small3.7bis}
\P\left(\sup_{(0,1)}\norm{Z}^2_2 \leq  {M} \right)\geq
C_p( {M} )\exp\left(-\frac{c_p"}{ {M} ^p} \right),
\end{equation}
Remark that for any $p$, $\eps_0$ we have $N(p, {M} )<\infty$.
Moreover, it is well-known that for any $d_0>0$, $P(d_0)>0$, which yields $C_p( {M} )>0$ and then
 Lemma \ref{PropSmallH3}.

\subsection{Proof of Lemma \ref{lem7.2}}\label{Preuve7.2}
\

\noindent For simplicity in the redaction, we restrict our
attention to the case $f=0$. The generalisation is easy.

Recall \eqref{Eq1.10bis}
$$
\E\abs{X_N(t)}^2\leq e^{- \mu_1 t}\abs{x_0}^2+\frac{c}{\mu_1}B_0.
$$
Since $(X^1,X^2)$ is a coupling of $(\Dr(X_N(\cdot,x_0^1)),\Dr(X_N(\cdot,x_0^1)))$ on $T\N$, we obtain
\begin{equation}\label{7.1}
\E\left(\abs{X^1(nT)}^2+\abs{X^2(nT)}^2\right)\leq e^{- \mu_1
nT}\left(\abs{x_0^1}^2+\abs{x_0^2}^2\right)+2\frac{c}{\mu_1}B_0.
\end{equation}
Since $(X^1,X^2)$ is a strong Markov process, it can be deduced
that there exist $C_6$ and
 $\gamma_6$ such that for any $x_0\in H$
\begin{equation}\label{7.2bis}
\E\left(e^{\gamma_6 \tau_{L^2}'}\right)\leq C_6\left(1+\abs{x_0^1}^2+\abs{x_0^2}^2\right),
\end{equation}
where
$$
\tau_{L^2}'=\inf\left\{t\in T\N\backslash\{0\}\,\left|\,   \abs{X^1(t)}^2+\abs{X^2(t)}^2\geq 4cB_0   \right.\right\}.
$$
Taking into account \eqref{7.2bis}, a standard argument gives
that, in order to establish Lemma \ref{lem7.2}, it is sufficient
to prove that there exist $(p_7,T_7)$ such that
\begin{equation}\label{7.2.1}
\P\left(\abs{X_N(t,x_0)}^2\leq \delta_3\right)\geq p_7(\delta_3,t)>0,
\end{equation}
provided $N\in\N$, $t\geq T_7(\delta_3)$ and $\abs{x_0}^2\leq 4cB_0$.

\noindent We set
$$
Z(t)=\int_0^te^{-A(t-s)}\phi(X_N(s))dW(s),\quad Y_N=X_N-P_NZ,\quad
M=\sup_{(0,t)}\norm{Z}_2^2.
$$

\noindent Assume that there exist
 $M_7(\delta_3)>0$ and $T_7(\delta_3)$ such that
\begin{equation}\label{7.2.2}
M\leq M_7(\delta_3)\quad\textrm{ implies }\quad
\abs{{Y_N}(t)}^2\leq \frac{\delta_3}{4},
\end{equation}
provided  $t\geq T_7(\delta_3)$ and $\abs{x_0}^2\leq 4cB_0$. Then
\eqref{7.2.1} results from Lemma \ref{PropSmallH3} with
$$
M=\min\left\{M_7(\delta_3),\frac{\delta_3}{4}\right\}.
$$
\noindent We now prove \eqref{7.2.2}. Remark that
\begin{equation}\label{7.2.3}
\frac{d}{dt}{Y_N}+  A {Y_N}+P_N B({Y_N}+{P_N Z})=0.
\end{equation}
Taking the scalar product of \eqref{7.2.3} with ${Y_N}$, it follows that
\begin{equation}\label{7.2.4bis}
\frac{d}{dt}\abs{{Y_N}}^2+2 \norm{{Y_N}}^2=-2({Y_N},B({Y_N}+{P_N Z})).
\end{equation}
Recalling that $(B(y,x),x)=0$, we obtain
$$
-2({Y_N},B({Y_N}+{P_N Z}))=-2({Y_N},({Y_N},\nabla){P_N Z})-2({Y_N},B({P_N Z})).
$$

\noindent We deduce from H\"older inequalities and Sobolev
embedding that
$$
-(z,(x,\nabla)y)\leq c\norm{z}\norm{x}\norm{y}.
$$
Hence it follows from \eqref{7.2.4bis} that
$$
\frac{d}{dt}\abs{{Y_N}}^2+2 \norm{{Y_N}}^2\leq c\norm{Z}^2\norm{{Y_N}}+c\norm{Z}\norm{{Y_N}}^2,
$$
which yields, by an arithmetico-geometric inequality,
$$
\frac{d}{dt}\abs{{Y_N}}^2+2 \norm{{Y_N}}^2\leq cM^\frac{1}{2}\norm{{Y_N}}^2+cM^\frac{3}{2}.
$$
It follows that $ M\leq \frac{1}{c^2}$ implies
\begin{equation}\label{7.2.4}
\frac{d}{dt}\abs{{Y_N}}^2+ \norm{{Y_N}}^2\leq cM^\frac{3}{2}\quad \textrm{ on }\; (0,t).
\end{equation}
Integrating, we deduce from $\abs{x_0}^2\leq 4cB_0$ that
$$
\abs{{Y_N}(t)}^2\leq 4c e^{-  \mu_1t}B_0+c\left(\frac{M^\frac{3}{2}}{  \mu_1}\right).
$$
Choosing $t$ sufficiently large and $M$ sufficiently small we
obtain \eqref{7.2.2} which yields \eqref{7.2.1} and then Lemma
\ref{lem7.2}.
\begin{Remark}
In order to avoid a lengthy proof, we have not splitted the
arguments in several cases as in the proof of Proposition
\ref{Prop_Lyap}. The reader can complete the details.
\end{Remark}

\subsection{Proof of Lemma \ref{lem7.4}}\label{Preuve7.4second}
\

 We use the decomposition $X_N=Y_N+P_NZ$ defined in section \ref{Preuve7.2} and set
$$
M=\sup_{(0,T)}\norm{Z}_2^2.
$$
 Integrating \eqref{7.2.4}, we obtain for $M$ satisfying the same
 assumption $M\leq\frac{1}{c^2}$
$$
\frac{1}{T}\int_0^T\norm{Y_N(t)}^2dt\leq \frac{1}{T}\abs{x_0}^2+cM^{\frac{3}{2}},
$$
which yields, by a Chebyshev inequality,
\begin{equation}\label{7.7.2}
\lambda\left(t\in (0,T)\,\left|\,\norm{Y_N(t)}^2\leq\frac{2}{T}\abs{x_0}^2+2cM^{\frac{3}{2}}\right.\right)\geq\frac{T}{2},
\end{equation}
where   $\lambda$ denotes the Lebesgue measure on $(0,T)$.

Setting
$$
\tau_{\H_1}=\inf\left\{t\in (0,T)\,\left|\, \norm{Y_N(t)}^2\leq
\frac{2}{T}\abs{x_0}^2+2cM^{\frac{3}{2}} \right.\right\},
$$
we deduce from \eqref{7.7.2} and the continuity of $Y_N$ that
\begin{equation}\label{7.7.3}
\norm{Y_N(\tau_{\H_1})}^2\leq
\frac{2}{T}\abs{x_0}^2+2cM^{\frac{3}{2}}.
\end{equation}

 Taking the scalar product of $2AY$ and \eqref{7.2.3}, we obtain
\begin{equation}\label{7.7.7}
\frac{d}{dt}\norm{Y_N}^2+2 \norm{Y_N}_2^2=-2(AY_N,B(Y_N+P_NZ)).
\end{equation}
It follows from H\"older inequalities, Sobolev Embeddings and Agmon inequality that
$$
-2(Ay,\widetilde B(x,z))\leq
c\norm{y}_2\norm{z}^\frac{1}{2}\norm{z}_2^\frac{1}{2}\norm{x},
$$
where $\widetilde B(x,y)=(x,\nabla)y+(y,\nabla)x$. Hence, we
obtain by applying arithmetico-geometric inequalities
$$
\Espace
\left\{
\begin{array}{lcccl}
-2(A{Y_N},B({Y_N}))&\leq& c\norm{{Y_N}}_2^\frac{3}{2}\norm{{Y_N}}^\frac{3}{2}&\leq &\frac{1}{4}\norm{{Y_N}}^2_2+c\norm{{Y_N}}^6,\\
-2(A{Y_N},B({P_N Z}))&\leq& c\norm {Y_N}_2\norm Z^\frac{3}{2}\norm Z_2^\frac{1}{2}&\leq&
\frac{1}{4}\norm{{Y_N}}^2_2+c\norm Z_2^4 ,\\
-2(A{Y_N},\widetilde B({Y_N},{P_N Z}))&\leq&  c\norm
{Y_N}_2^\frac{3}{2}\norm {Y_N}^\frac{1}{2}\norm Z&\leq &
c\norm Z\norm {Y_N}_2^2
\end{array}
\right.
$$
Remarking that $B({Y_N}+{P_N Z})=B({Y_N})+\widetilde B({Y_N},{P_N
Z})+B({P_N Z})$, it follows from \eqref{7.7.7}  that
$M\leq\frac{1}{4c}$ implies
\begin{equation}\label{7.7.8}
\frac{d}{dt}\norm{{Y_N}}^2+ \norm{{Y_N}}_2^2\leq
c\norm{{Y_N}}^2\left(\norm{{Y_N}}^4-4K_0^2\right) +cM^2,
\end{equation}
where $K_0$ is defined in \eqref{2.9bis}. Let us set
$$
\sig_{\H_1}=\inf\left\{t\in(\tau_{\H_1},T)\,\left|\,  \norm{{Y_N}(t)}^2> 2 K_0   \right.\right\},
$$
and remark that on $(\tau_{\H_1},\sig_{\H_1})$, we have
\begin{equation}\label{7.7.9}
\frac{d}{dt}\norm{{Y_N}}^2+ \norm{{Y_N}}_2^2\leq
cM^2.
\end{equation}
Integrating, we obtain that
\begin{equation}\label{7.7.10}
\norm{{Y_N}(\sig_{\H_1})}^2+ \int_{\tau_{\H_1}}^{\sig_{\H_1}}\norm{{Y_N}(t)}_2^2dt\leq \norm{{Y_N}(\tau_{\H_1})}^2
+cM^2.
\end{equation}
Combining \eqref{7.7.3} and \eqref{7.7.10},  we obtain that, for $M$ and $\abs{x_0}^2$ sufficiently small,
$$
\norm{{Y_N}(\sig_{\H_1})}^2\leq \frac{\delta_4}{4}\wedge K_0,
$$
which yields $\sig_{\H_1}= T$. It follows that
\begin{equation}\label{7.7.6}
\norm{ X_N( T)}^2\leq \delta_4,
\end{equation}
provided $M$ and $\abs{x_0}^2$ sufficiently small. It remains to
use Lemma \ref{PropSmallH3} to get Lemma \ref{lem7.4}.

\subsection{Proof of Lemma \ref{lem7.5}}\label{Preuve7.5}
\

  It follows from \eqref{7.7.9} that
$$
 \int_0^T\norm{Y_N(t)}_2^2dt\leq \norm{x_0}^2+
cM^2,
$$
provided $M\leq \frac{1}{4c}$ and
$
\norm{x_0}^2+
cM^2\leq K_0.
$

\noindent   Applying the same argument as in the previous
subsection, it is easy to deduce that
 there exists a stopping times $\tau_{\H_2}\in (0,T)$ such that
\begin{equation}\label{7.7.3bis}
\norm{Y_N(\tau_{\H_2})}_2^2\leq \frac{2}{T}\left(\norm{x_0}^2+
cM^2\right),
\end{equation}
provided $M$ and $\norm{x_0}$ are sufficiently small.

 Taking the scalar product of \eqref{7.2.3} and $2A^2{Y_N}$, we obtain
\begin{equation}\label{7.15}
\frac{d}{dt}\norm{{Y_N}}_2^2+2 \norm{{Y_N}}_3^2= -2\left(A^\frac{3}{2}{Y_N},A^\frac{1}{2}B({Y_N}+{P_N Z})\right).
\end{equation}
Applying H\"older inequality, Sobolev Embeddings $\H_2\subset L^\infty$ and $\H_1\subset L^4$
 and arithmetico-geometric inequality, we obtain
$$
-2\left(A^\frac{3}{2}y,A^\frac{1}{2}B(x,y)\right)\leq
c\norm{y}_3\norm{x}_2\norm{y}_2\leq
 \frac{1}{4} \norm{y}_3^2+c\left(\norm{x}_2^4+\norm{y}_2^4\right).
$$
Hence we deduce  from \eqref{7.15} and from $B({Y_N}+{P_N Z})=B({Y_N})+B({Y_N},{P_N Z})+B({P_N Z})$
\begin{equation}\label{7.16}
\frac{d}{dt}\norm{{Y_N}}_2^2+ \norm{{Y_N}}_3^2\leq
c\norm{{Y_N}}_2^2(\norm{{Y_N}}_2^2-2K_1)+c\norm{Z}_2^4 ,
\end{equation}
where $K_1$ is defined as $K_0$ in \eqref{2.9bis} but with a
different $c$. We set
$$
\sig_{\H_2}=\inf\left\{t\in(\tau_{\H_2},T)\,\left|\,  \norm{{Y_N}(t)}_2^2> 2 K_1   \right.\right\},
$$
Integrating \eqref{7.16}, we obtain
$$
\norm{{Y_N}(\sig_{\H_2})}_2^2+ \int_{\tau_{\H_2}}^{\sig_{\H_2}}\norm{{Y_N}(t)}_3^2dt\leq
\norm{{Y_N}(\tau_{\H_2})}_2^2+cM^2.
$$
Taking into account \eqref{7.7.3bis} and choosing $\norm{x_0}^2$ and $M^2$ sufficiently small,
 we obtain
$$
\norm{{Y_N}(\sig_{\H_2})}_2^2\leq \frac{\delta}{4}\wedge K_1.
$$
It follows that $\sig_{\H_2}=T$ and that
\begin{equation}\label{7.7.6bis}
\norm{ X_N( T)}^2\leq \delta,
\end{equation}
provided $M$ and $\norm{x_0}$ sufficiently small, which yields \eqref{2.5}.

\section{Proof of Theorem \ref{Th}}

As already explained, Theorem \ref{Th} follows from Proposition
\ref{ThN}.
 We now prove Proposition \ref{ThN}. Let $(x_0^1,x_0^2)\in(\H_2)^2$.
Let us recall that the process $(X^1,X^2)$ is defined at the beginning of section $3$.

Let $\delta>0$, $T\in (0,1)$ be as in Proposition \ref{Prop_I_0}
and $\tau$ defined in \eqref{3.1bis}, setting
$$
\tau_1=\tau,\quad \tau_{k+1}=\inf\left\{t>\tau_{k}\,\left|\,   \norm{X^1(t)}_2^2\vee\norm{X^2(t)}_2^2\leq \delta  \right.\right\}.
$$
it can be deduced from the strong Markov property of $(X^1,X^2)$
and from Proposition \ref{Prop_Lyap} that
$$
\E\left(e^{\alpha \tau_{k+1}}\right)\leq K"\E\left(e^{\alpha \tau_k}\left(1+\abs{X^1(\tau_k)}^2+\abs{X^2(\tau_k)}^2\right)\right),
$$
which yields, by the Poincar\'e inequality,
$$
\Espace
\left\{
\begin{array}{rcl}
\E\left(e^{\alpha \tau_{k+1}}\right)&\leq& cK"(1+2\delta)\E\left(e^{\alpha \tau_k}\right),\\
\E\left(e^{\alpha \tau_{1}}\right)&\leq& K"\left(1+\abs{x_0^1}^2+\abs{x_0^2}^2\right).
\end{array}
\right.
$$
It follows that there exists $K>0$ such that
$$
\E\left(e^{\alpha \tau_{k}}\right)\leq K^k\left(1+\abs{x_0^1}^2+\abs{x_0^2}^2\right).
$$
Hence, applying Jensen inequality, we obtain that, for any $\theta\in(0,1)$
\begin{equation}\label{tau_k}
\E\left(e^{\theta\alpha \tau_{k}}\right)\leq K^{\theta k}\left(1+\abs{x_0^1}^2+\abs{x_0^2}^2\right).
\end{equation}
We deduce from Proposition \ref{Prop_I_0} and from \eqref{toy5bisbis} that
$$
\P\left(X^1(T)\not=X^2(T)\right)\leq\frac{1}{4},
$$
provided $(x_0^1,x_0^2)$ are in the ball of $(\H^2)^2$ with radius
$\delta$.

Setting
$$
k_0=\inf\left\{k\in\N\,\left|\,X^1(\tau_k+T)=X^2(\tau_k+T)\right.\right\},
$$
it follows that $k_0<\infty$ almost surely and that
\begin{equation}\label{k_0}
\P\left(k_0>n\right)\leq \left(\frac14\right)^n.
\end{equation}
Let $\theta\in (0,1)$. We deduce from Schwartz inequality that
$$
\E\left(e^{\frac{\theta}{2}\alpha \tau_{k_0}}\right)
=\sum_{n=1}^\infty\E\left(e^{\frac{\theta}{2}
\alpha\tau_{n}}1_{k_0=n}\right) \leq\sum_{n=1}^\infty
\sqrt{\P\left(k_0\geq n\right)\E\left(e^{\theta\alpha
\tau_{n}}\right)}.
$$
Combining \eqref{tau_k} and \eqref{k_0}, we deduce
$$
\E\left(e^{\frac{\theta}{2}\alpha \tau_{k_0}}\right)\leq
 \left(\sum_{n=0}^\infty\left(\frac{K^\theta}{2}\right)^n\right)
\left(1+\abs{x_0^1}^2+\abs{x_0^2}^2\right).
$$
Hence, choosing $\theta\in(0,1)$ sufficiently small, we obtain that there exists $\gamma>0$ non depending on $N\in\N$
such that
\begin{equation}\label{exp}
\E\left(e^{\gamma\tau_{k_0}}\right)\leq 4\left(1+\abs{x_0^1}^2+\abs{x_0^2}^2\right).
\end{equation}
Recall that if $(X^1,X^2)$ are coupled at time $t\in T\N$, then
they remain coupled for any time after. Hence $X^1(t)=X^2(t)$ for $t>\tau_{k_0}$.
It follows
$$
\P\left(X^1(nT)\not=X^2(nT)\right)\leq 4e^{-\gamma nT}\left(1+\abs{x_0^1}^2+\abs{x_0^2}^2\right).
$$
 Since $(X^1(nT),X^2(nT))$ is a coupling of
$((\Pcal_{nT}^N)^*\delta_{x_0^1},(\Pcal_{nT}^N)^*\delta_{x_0^2})$, we deduce from Lemma \ref{lem_coupling}
 \begin{equation}\label{2.1}
\norm{\left(\Pcal_{nT}^N\right)^*\delta_{x_0^1}-\left(\Pcal_{nT}^N\right)^*\delta_{x_0^2}}_{var}
\leq 4 e^{-\gamma n T}\left(1+\abs{x_0^1}^2+\abs{x_0^2}^2\right),
\end{equation}
for any $n\in\N$ and any $(x_0^1,x_0^2)\in(\H_2)^2$.

Recall that the existence of an invariant measure $\mu_N\in P(P_N H)$ is justified in section 1.3.
 Let $\lambda\in P(H)$ and $t\in\R^+$. We set $n=\lfloor\frac{t}{T}\rfloor$ and $C=4e^{\gamma T}$.
 Integrating $(x_0^1,x_0^2)$ over $((\Pcal_{t-nT}^N)^*\lambda)\otimes\mu_N$ in \eqref{2.1}, we obtain
$$
\norm{\left(\Pcal_{t}^N\right)^*\lambda-\mu_N}_{var}
\leq C e^{-\gamma t}\left(1+\int_{H}\abs{x_0}^2\,\lambda(dx)\right),
$$
which establishes \eqref{EqThN}.

\appendix

\section{ Proof of \eqref{2.15}}\label{Preuve_2.15}

\noindent For simplicity in the redaction, we omit $\theta$ and $N$ in our notations.

\noindent Remark that
\begin{equation}\label{2.12}
J=\left(\nabla \E\left(g(X(T))\psi_X\right),h\right)=J_1+2J_2,
\end{equation}
where
$$
\Espace
\left\{
\begin{array}{lcl}
J_1&=&\E\left(\left(\nabla g(X(T)),\eta(T,0)\cdot h\right)\psi_X\right),\\
J_2&=&\E\left( g(X(T))\psi'_X\int_0^T\left(AX(t),A(\eta(t,0)\cdot h)\right)ds\right).
\end{array}
\right.
$$
According to \cite{Nualart}, let us denote by $D_s F$  the
Malliavin derivative of $F$ at time $s$. We have the following
formula of the Malliavin derivative of the solution of a
stochastic differential equation
$$
 D_sX(t)=1_{t\geq
s}\eta(t,s)\cdot\phi(X(s)),
$$
which yields
\begin{equation}\label{A2}
\int_0^t D_sX(t)\cdot m(s)\,ds=G(t)\cdot m,
\end{equation}
where
$$
G(t) \cdot m= \int_0^t\eta(t,s)\cdot\phi(X(s))\cdot m(s)\,ds.
$$
The uniqueness of the solutions gives
$$
\eta(t,0)\cdot h=\eta(t,s)\cdot (\eta(s,0)\cdot h) \textrm{ for any } 0\leq s \leq t,
$$
which yields
$$
\eta(T,0)\cdot h=\frac{1}{T}\int_0^T\eta(t,s)\cdot (\eta(s,0)\cdot
h)\,ds.
$$
 Setting
$$
w(s)=\phi^{-1}(X(s))\cdot \eta(s,0)\cdot h,
$$
we infer from \eqref{A2}
\begin{equation}\label{A3}
\eta(T,0)\cdot h=\frac{1}{T}\,G(T)\cdot
w=\frac{1}{T}\int_0^TD_sX(T)\cdot w\,ds,
\end{equation}
which yields
$$
\left(\nabla g(X(T)),\eta(T,0)\cdot
h\right)
=\frac{1}{T}\int_0^T\left(\nabla g(X(T)),D_sX(T)\cdot
w\right)ds.
$$
Remark that
$$
\left(D_sg(X(T)),w\right)=\left(\nabla g(X(T)),D_sX(T)\cdot
w\right).
$$
It follows
$$
\left(\nabla g(X(T)),\eta(T,0)\cdot
h\right)=\frac{1}{T}\int_0^T\left(D_sg(X(T)),w\right)ds,
$$
 which yields
\begin{equation}\label{2.13ter}
J_1=\frac{1}{T}\E\int_0^T \psi_X\left(D_{s} g(X(T)), w\right)ds.
\end{equation}
Recall that the Skohorod integral is the dual operator of the
Malliavin derivative (See \cite{Nualart}). It follows
\begin{equation}\label{2.14}
J_1=\frac{1}{T}\E \left(g(X(T))\int_0^T \psi_X(w(t),dW(t))\right).
\end{equation}
Recall the formula of integration of a product
\begin{equation}\label{6.1}
\int_0^T \psi_X(w(t),dW(t))=\psi_X \int_0^T (w(t),dW(t))- \int_0^T
\left(D_{s} \psi_X,w(s)\right)ds.
\end{equation}
Remark that
$$
D_s\psi_X=2\psi_X'\int_0^T AD_sX(t) \cdot (AX(t))\,dt,
$$
which yields, by $\left(AX(t),A D_{s}X(t)\cdot
w(s)\right)=\left(w(s),A D_{s}X(t)\cdot (AX(t))\right)$,
$$
 \int_0^T\left(D_{s} \psi_X,w(s)\right)ds=2\psi'_X \int_0^T\int_0^T \left(AX(t),A D_{s}X(t)\cdot w(s)\right)\,dtds.
$$
 We deduce from \eqref{A3} that
\begin{equation}\label{6.2}
 \int_0^T\left(D_{s} \psi_X,w(s)\right)ds=2\psi'_X \int_0^T t\left(AX(t),A \eta(t,0)\cdot h\right)\,dt.
\end{equation}
Remark that
$$
\psi_X=\psi'_X=0\quad\textrm{ if }\quad\sig<T.
$$
Hence combining   \eqref{6.1} and \eqref{6.2}, we obtain
$$
\int_0^T \psi_X(w(t),dW(t))=\psi_X \int_0^\sig (w(t),dW(t))-
2\psi'_X \int_0^\sig t\left(AX(t),A \eta(t,0)\cdot h\right)\,dt.
$$
Thus, \eqref{2.15} follows from \eqref{2.12} and  \eqref{2.14}.



\BLANC{


}

\end{document}